\documentclass[11pt]{article}
\pdfoutput=1
\usepackage{amssymb,amsfonts,amsmath,amsthm,amscd,dsfont,mathrsfs}
\usepackage{graphicx,float,psfrag,epsfig,color}
\usepackage{latexsym}

\def\E{{\mathbb E}}

\def\onu{\overline{\nu}}
\def\E{\mathbb{E}}

\def\de{{\rm d}}

\newcommand{\targ}{\xo}
\newcommand{\bhat}{\hat{x}}
\newcommand{\bhatee}{\bhat^t}
\newcommand{\bhatpo}{\bhat^{t+1}}

\def\de{{\rm d}}

\def\reals{{\mathbb R}}
\def\E{{\mathbb E}}
\def\<{\langle}
\def\>{\rangle}

\def\prob{{\mathbb P}}

\newcommand{\beqa}{\begin{eqnarray}}
\newcommand{\eeqa}{\end{eqnarray}}
\newcommand{\rhoMSE}{\rho_{\mbox{\rm \tiny MSE}}}

\newcommand{\bitem}{\begin{itemize}}
\newcommand{\eitem}{\end{itemize}}
\newcommand{\beq}{\begin{equation}}
\newcommand{\eeq}{\end{equation}}
\newcommand{\goto}{\rightarrow}

\newcommand{\cE}{{\cal E}}

\newcommand{\eps}{{\varepsilon}}

\newcommand{\BPDN}{P_{2,\lambda,1}}
\newcommand{\LASSO}{LASSO}
\newcommand{\BP}{P_{1}}
\newcommand{\hash}{\pm}

\newcommand{\lslim}{\mbox{\rm ls lim}}

\definecolor{cardinal}{rgb}{.64,0.,.11}

\definecolor{Red}{rgb}{1,0,0}
\definecolor{Blue}{rgb}{0,0,1}
\definecolor{Olive}{rgb}{0.41,0.55,0.13}
\definecolor{Green}{rgb}{0,1,0}
\definecolor{MGreen}{rgb}{0,0.8,0}
\definecolor{DGreen}{rgb}{0,0.55,0}
\definecolor{Yellow}{rgb}{1,1,0}
\definecolor{Cyan}{rgb}{0,1,1}
\definecolor{Magenta}{rgb}{1,0,1}
\definecolor{Orange}{rgb}{1,.5,0}
\definecolor{Violet}{rgb}{.5,0,.5}
\definecolor{Purple}{rgb}{.75,0,.25}
\definecolor{Brown}{rgb}{.75,.5,.25}
\definecolor{Grey}{rgb}{.5,.5,.5}
\definecolor{Pink}{rgb}{1,0,1}
\definecolor{DBrown}{rgb}{.5,.34,.16}

\definecolor{Black}{rgb}{0,0,0}

\numberwithin{equation}{section}
\newtheorem{theorem}{Theorem}[section]

\newtheorem{lemma}{Lemma}[section]
\newtheorem{proposition}{Proposition}[section]
\newtheorem{definition}{Definition}[section]
\newtheorem{finding}[theorem]{Finding}

\newcommand{\cF}{{\cal F}}

\def\normal{{\sf N}}
\def\E{{\mathbb E}}
\def\prob{{\mathbb P}}
\def\reals{{\mathbb R}}

\def\MSE{{\rm MSE}}
\def\NPI{{\rm NPI}}
\def\LSF{{\rm LSF}}
\def\FMSE{{\rm fMSE}}
\def\FNPI{{\rm fNPI}}

\def\cJ{{\cal J}}

\def\Ms{M}

\def\HFP{{\rm HFP}}
\def\SC{{\rm SC}}

\def\FDR{{\rm FDR}}
\def\FDeR{{\rm FDeR}}
\def\DR{{\rm DR}}

\def\MSE{{\rm MSE}}
\def\eMSE{{\rm eMSE}}
\def\MSR{{\rm MSR}}
\def\MAE{{\rm MAE}}

\def\EqMSE{{\rm EqMSE}}
\def\EqMSR{{\rm EqMSR}}
\def\EqMAE{{\rm EqMAE}}
\def\EqDR{{\rm EqDR}}
\def\EqPMSR{{\rm EqPMSR}}

\def\AMPT{{\rm AMPT}}
\def\LASSO{{\rm LASSO}}

\def\onu{\overline{\nu}}
\def\df{{\rm df}}

\newcommand{\stMSE}{{\sf mse}}
\newcommand{\npi}{{\sf npi}}
\newcommand{\Formal}{{\sc\rm Formal}}

\title{The Noise-Sensitivity Phase Transition in Compressed Sensing}

\author{David L. Donoho\thanks{Department of Statistics, Stanford University},
\;\;\;
Arian Maleki\thanks{Department of Electrical Engineering, Stanford University},
\;\; and \;\; Andrea~Montanari${}^{*,\dagger}$}

\newcommand{\xo}{x^0}
\newcommand{\hxl}{\hat{x}^{1,\lambda}}

\begin{document}

\maketitle

\begin{abstract}
Consider the noisy underdetermined system
of linear equations: $y=Ax^0 + z^0$, with $n \times N$
measurement matrix $A$, $n < N$, and Gaussian white noise
$z^0 \sim \normal(0,\sigma^2 I )$.
Both $y$ and $A$ are known, both $x^0$ and $z^0$ are unknown, and
we seek an approximation to $x^0$.

When $x^0$ has
few nonzeros, useful approximations are often obtained by $\ell_1$-penalized
$\ell_2$ minimization, in which the reconstruction $\hxl$
solves $ \min \| y - Ax\|_2^2/2 +  \lambda \|x\|_1$.

 Evaluate performance by mean-squared error ($\MSE =\E ||\hxl - \xo||_2^2/N$).
Consider matrices $A$ with iid Gaussian entries and a large-system limit in which
$n,N\to\infty$ with $n/N \to  \delta$ and $k/n \to \rho$.
Call   the ratio $\MSE/\sigma^2$ the {\it noise sensitivity}.
We develop formal expressions for the MSE
 of $\hxl$,
and evaluate its worst-case formal noise sensitivity  over
all types of $k$-sparse signals.
The phase space $ 0 \leq \delta, \rho \leq 1$
is partitioned by curve $\rho = \rhoMSE(\delta)$
into two regions. Formal noise sensitivity
is bounded throughout the region $\rho < \rhoMSE(\delta)$ and is
unbounded throughout the region $\rho > \rhoMSE(\delta)$.

The phase boundary $\rho = \rhoMSE(\delta)$
is {\it identical to} the previously-known phase transition curve
for equivalence of  $\ell_1 - \ell_0$ minimization in the
$k$-sparse noiseless case. Hence
a single phase boundary describes the fundamental phase transitions both for the
noiseless and noisy cases.

 Extensive computational experiments
validate the predictions of this formalism, including the existence
of game theoretical structures underlying it (saddlepoints in the payoff,
least-favorable signals and maximin penalization).

Underlying our  formalism is an
approximate message passing soft thresholding algorithm (AMP)
introduced earlier by the authors. Other papers by the authors detail expressions
for the formal MSE of AMP and its close connection to
$\ell_1$-penalized  reconstruction.
Here we derive the minimax formal MSE of AMP
and then read out results for $\ell_1$-penalized
reconstruction.
\end{abstract}

\vspace{.1in}
Lasso. Basis Pursuit. Minimax Risk over Nearly-Black Objects.
Minimax Risk of Soft Thresholding.
\vspace{.1in}

\vspace{.1in}
{\bf Acknowledgements.}
Work partially supported by
NSF DMS-0505303, NSF DMS-0806211, NSF CAREER CCF-0743978.
Thanks to Iain Johnstone and Jared Tanner for helpful discussions.

\vspace{.1in}

\newpage


\section{Introduction}

Consider the noisy underdetermined system
of linear equations:
\beq
\label{eq:obsdata} y=Ax^0 + z^0\, ,
\eeq
where the matrix $A$ is $n\times N$, $n < N$,
the $N$-vector $x^0$ is $k$-sparse (i.e. it has at most $k$
non-zero entries), and $z^0 \in \reals^{n}$ is
a Gaussian white noise $z^0\sim \normal(0,\sigma^2 I )$.
Both $y$ and $A$ are known, both $x^0$ and $z^0$ are unknown, and
we seek an approximation to $x^0$.

A very popular approach estimates $x^0$ via the solution $x^{1,\lambda}$
of the following convex optimization problem
\beq \label{BPDN}
(\BPDN) \qquad  {\rm minimize} \quad  \frac{1}{2}\, \| y - Ax\|_2^2 +
\lambda \|x\|_1 .
 \eeq
Thousands of articles use or study this approach,
which has variously been called LASSO, Basis Pursuit, or more
prosaically, $\ell_1$-penalized least-squares \cite{Tibs96,BP95,BP}.
 There is a clear need to understand the extent to which $(\BPDN)$
 accurately recovers $\xo$.   Dozens of papers
 present partial results, setting forth often loose
 bounds on the behavior of $\hxl$ (more below).

Even in the noiseless case $z^0 = 0$, understanding the reconstruction problem
 (\ref{eq:obsdata}) poses a challenge,
as  the underlying system of equations $y = Ax^0$ is
 underdetermined.
In  this case
it is informative to consider  $\ell_1$ minimization,
 \begin{align} \label{BP}
(\BP) \qquad &\mbox{minimize } \|x\|_1\, ,\\
& \mbox{ subject to } y =  Ax .
 \end{align}
\newcommand{\hxz}{\hat{x}^{1,0}}
This is the $\lambda=0$ limit of (\ref{BPDN}):
its solution obeys $\hxz=\lim_{\lambda\to 0}x^{1,\lambda}$.

The most precise information about behavior of  $\hxz$
is obtained by large-system analysis;
let $n,N$ tend to infinity so that\footnote{Here and below
we write $a\sim b$ if $a/b\to 1$ as both quantities tend to
infinity.} $n \sim \delta N$
and correspondingly let the number of nonzeros $k \sim \rho n$;
thus we have a phase space $0 \leq \delta,\rho \leq 1$, expressing
different combinations of undersampling $\delta$ and sparsity $\rho$.
When the matrix $A$ has iid Gaussian elements,
phase space $ 0 \leq \delta,\rho \leq 1$ can be divided into
two components, or {\it phases}, separated by a curve $\rho = \rho_{\ell_1}(\delta)$,
which can be explicitly computed.
Below this curve, $x^0$ is sufficiently sparse that
$\hxz = x^0$ with high probability
and therefore $\ell_1$ minimization perfectly recovers the sparse vector $\hxz$.
Above this curve, sparsity is not sufficient:
we have $\hxz \neq x^0$ with high probability.
Hence the curve $\rho = \rho_{\ell_1}(\delta)$, $0 < \delta < 1$,
indicates the precise tradeoff  between undersampling
and sparsity.

Many authors have considered the behavior of $\hxl$ in the noisy case
but results are somewhat less conclusive.
The most well-known analytic approach is the Restricted Isometry Principle (RIP),
developed by Cand\`es and Tao \cite{CandesTao,CandesTaoDantzig}.
Again in the case where $A$ has iid Gaussian entries, and in the same large-system limit,
the RIP implies that, under sufficient sparsity of $x^0$, with high probability
one has stability bounds of the form $\| \hxl - x^0 \|_2 \leq
C(\delta,\rho) \|z^0\|_2\, \log N$.
The region where $C(\delta,\rho) < \infty $ was orginally an implicitly known,
but clearly nonempty region of the $(\delta,\rho)$
phase space.
Blanchard, Cartis and Tanner \cite{BlCaTa09}
recently improved the estimates of $C$ in the case of Gaussian matrices $A$,
by careful large deviations analysis,
and by developing an asymmetric RIP, obtaining the largest
region where $\hxl$ is currently known to be stable.
Unfortunately  as they show, this region is still relatively small compared to
the region $\rho < \rho_{\ell_1}(\delta)$, $0 < \delta < 1$.

It may seem that, in the presence of noise,
the precise tradeoff between undersampling and sparsity
worsens dramatically, compared to the noiseless case.
In fact, the opposite is true.  In this paper, we show that in the presence
of Gaussian white noise, the mean-squared error of the optimally tuned
$\ell_1$ penalized least squares estimator behaves well
over quite a large region of the phase plane, in fact,
it is finite over the exact same region of the phase plane
as the region of $\ell_1-\ell_0$ equivalence derived in the noiseless case.

Our main results, stated in Section \ref{sec:MainResults},
give explicit evaluations for the the worst-case
formal mean square error of $\hxl$ under given conditions of noise,
sparsity and undersampling.  Our results indicate the
noise sensitivity of solutions to (\ref{BPDN}), the optimal
penalization parameter $\lambda$,
and the hardest-to-recover sparse vector.
As we show, the noise sensitivity
exhibits a phase transition in the undersampling-sparsity $(\delta,\rho)$
domain along a curve $\rho = \rhoMSE(\delta)$, and this curve
is precisely the same as the $\ell_1$-$\ell_0$ equivalence curve
$\rho_{\ell_1}$.

Our results might be compared to work of Xu and Hassibi
\cite{XuHassibi}, who considered
a different departure from the noiseless case. In their work, the noise $z^0$ was
still vanishing, but the vector $x_0$ was allowed to be an $\ell_1$-norm bounded perturbation to
a $k$-sparse vector.  They considered stable recovery with respect to such small perturbations
and showed that the natural boundary for such stable recovery is again the curve $\rho = \rhoMSE(\delta)$.

\subsection{Results of our Formalism}

We define below a so-called formal MSE ($\FMSE$), and evaluate the (minimax, formal) {\it noise sensitivity}:
\begin{equation} \label{eq:DefMStar}
     M^*(\delta,\rho) = \sup_{\sigma > 0}  \max_\nu    \min_\lambda \FMSE(\hxl,\nu,\sigma^2)/\sigma^2;
\end{equation}
here $\nu$ denotes the marginal distribution of $x^0$ (which has
fraction of nonzeros not larger than $\rho  \delta$),
and $\lambda$ denotes the tuning parameter of the $\ell_1$-penalized $\ell_2$
minimization.
Let $M^\hash(\eps)$ denote the minimax MSE of scalar thresholding,
defined in Section 2 below.  Let $\rhoMSE(\delta)$ denote the solution of
 \begin{eqnarray} \label{DefRhoMSE}
M^\hash(\rho\delta) = \delta\, .
\end{eqnarray}
Our main {\it theoretical} result is the formula
\beq \label{eqMStarResult}
   M^*(\delta,\rho) =  \left\{ \begin{array}{ll}
                                         \frac{M^\hash(\delta\rho)}{1-M^\hash(\delta\rho)/\delta} , & \rho < \rhoMSE(\delta) ,\\
                                         \infty  ,& \rho \geq  \rhoMSE(\delta). \\
                                           \end{array} \right.
\eeq
\begin{figure}
\center{\includegraphics[width=10.cm]{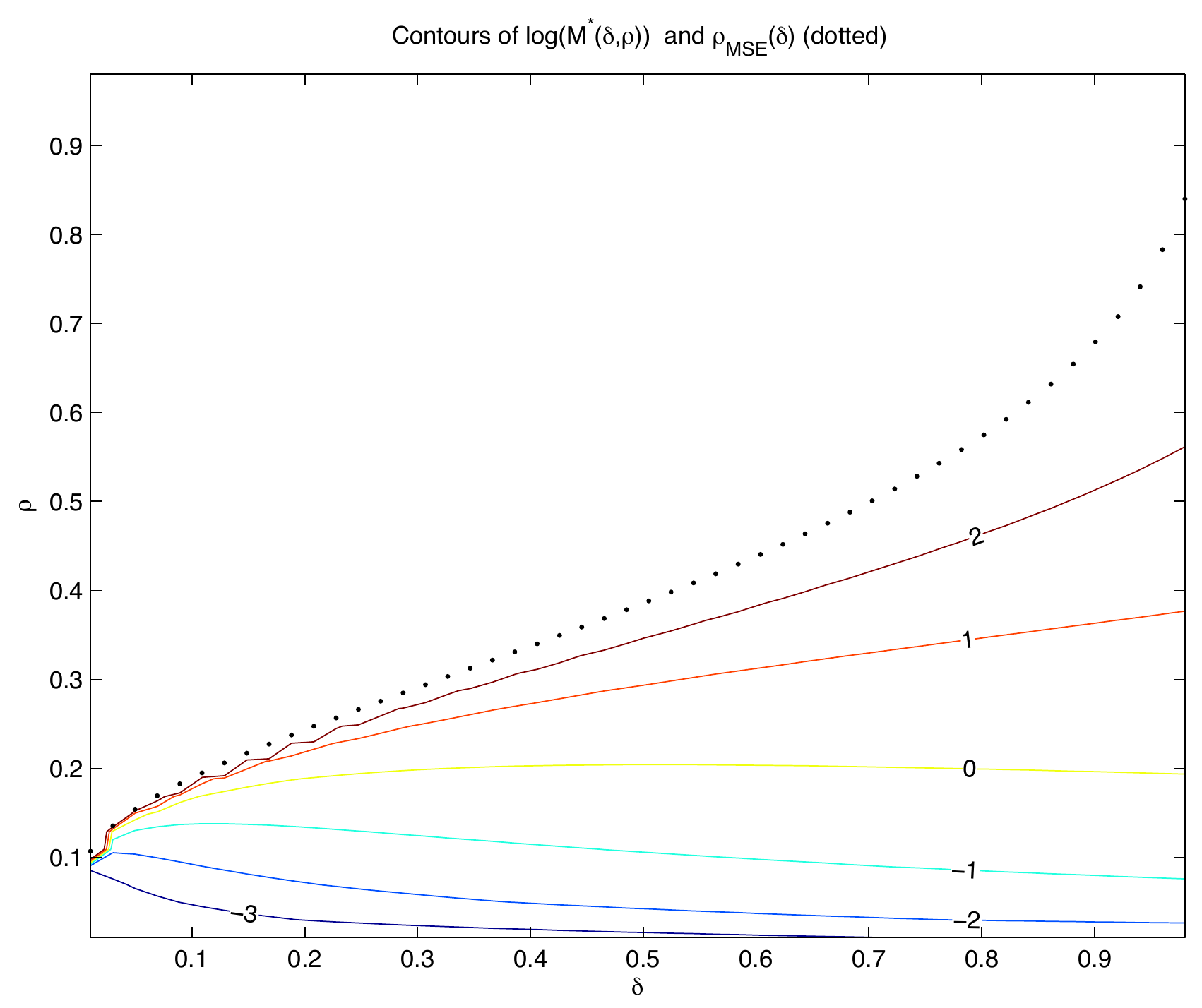}}
\caption{{\small Contour lines of the minimax noise sensitivity
$\Ms^*(\delta, \rho)$ in the $(\rho,\delta)$ plane. The
dotted black curve graphs the phase boundary $(\delta,\rhoMSE(\delta))$.
Above this curve, $\Ms^*(\delta, \rho) =\infty$. The
colored lines present level sets of $\Ms^*(\delta, \rho)
=1/8$, $1/4$, $1/2$, $1$, $2$, $4$ (from bottom to top).}}
\label{fig:Figure001}
\end{figure}

Quantity (\ref{eq:DefMStar})
is the payoff of a traditional two-person zero sum game, in which the undersampling and sparsity
are fixed in advance, the researcher
plays against Nature, Nature picks both a noise level and a signal distribution,
and the researcher picks a penalization level, in knowledge of Nature's
choices.  It is traditional in analyzing such games to identify the least-favorable
strategy of Nature (who maximizes payout from the researcher),
and the optimal strategy for the researcher (who wants to minimize payout).
We are able to identify both and give explicit formulas for the so-called
saddlepoint strategy, where Nature plays the least-favorable strategy
against the researcher and the researcher minimizes the consequent damage.  In Proposition \ref{prop:main}
below we give formulas for this pair of strategies.  The phase-transition structure evident in
(\ref{eqMStarResult}) is saying that above the curve $\rhoMSE$, Nature has available
unboundedly good strategies, to which the researcher has no effective response.

\subsection{Structure of the Formalism}

 Our approach is presented in Section \ref{sec:Formalism},
and uses a combination of ideas from decision theory in mathematical
statistics,
and message passing algorithms in information theory.
On the one hand, as already evident from formula (\ref{eqMStarResult}),
quantities from mathematical statistics play a key role in our formulas.
But since these quantities concern a completely different estimator in
a completely different problem -- the behavior of soft thresholding
in estimating a single normal mean, likely to be zero -- the superficial
appearance of the formulas conceals the type of analysis we are doing.
That analysis concerns the properties of an iterative soft thresholding
scheme introduced by the authors in \cite{DMM}, and further developed
here.  Our formalism neatly describes properties of the formal MSE
of AMP as expectations taken in the equilibrium states of a state evolution.
As described in \cite{TheoPredLassoOpChar}, we can calibrate AMP to have
the same operating characteristics as $\ell_1$-penalized least squares,
and by recalibration of the minimax formal MSE for AMP, we get
the above results.

\subsection{Empirical Validation}

We use the word {\it formalism} for the machinery underlying
our derivations because it is not (yet) a rigorously-proven method
which is known to give correct results under established regularity conditions.
In this sense our method has similarities to the replica and cavity methods
of statistical physics, famously useful tools without rigorous general
justification.

Our theoretical results are
validated here by computational experiments which show that
the predictions of our formulas are accurate, and, even more
importantly, that the underlying formal structure leading to
our predictions -- least-favorable objects, game-theoretic saddlepoints of the
MSE payoff function, maximin tuning of $\lambda$, unboundedness
of the noise sensitivity above phase transition-- can all be
observed experimentally.  Because our formalism makes
so many different kinds of predictions about quantities with clear
operational significance and about their dynamical evolution
in the AMP algorithm, it is quite different than some other formalisms,
such as the replica method, in which many fewer checkable predictions are made.
In particular, as demonstrated in \cite{DMM}, the present
formalism describes precisely the evolution of an actual low complexity algorithm.

Admittedly, by computational
means we can only check individual predictions in specific cases,
whereas a full proof could cover all such cases. However, we make available software
which checks these features so that interested researchers
can check the same phenomena at parameter values that we did
not investigate here.  {\it The evidence of our simulations
is strong; it is not a realistic possibility
that $\ell^1$-penalized least squares fails to have the
limit behavior discovered here.}

We focused in this paper on measurement matrices $A$ with Gaussian iid entries.
It was recently proved that the state evolution formalism
 at the core of our analysis is indeed asymptotically
correct for Gaussian matrices $A$ \cite{BayatiMontanari}.
We believe that similar results hold for matrices $A$ with uniformly bounded iid entries with zero
mean and variance $1/n$. However our results should extend to a broader
universality class including matrices with iid entries with same mean and
variance, under an appropriate light tail condition.
It is an outstanding mathematical challenge to prove
that such predictions are indeed correct for
a broader universality class of estimation problems.

As discussed in Section \ref{sec:StatPhys},
an alternative route also from statistical physics,
using the replica method has been recently used
to investigate similar questions.
We will argue that the present framework which makes predictions
about actual dynamical behavior of algorithms,
is computationally verifiable in great detail,
whereas the replica method itself applies to no constructive algorithm and makes
comparatively many fewer predictions.

%
%
\section{Minimax MSE of Soft Thresholding}
\label{sec:Scalar}

We briefly recall notions from, e.g.,  \cite{DJHS92,DJ94a} and then
generalize them.
We wish to recover an $N$ vector $x^0 = (x^0(i): 1 \leq i \leq N )$
which is observed in Gaussian white noise
\[
     y(i) = x^0(i) +  z^0(i), \qquad  1 \leq i \leq N,
\]
with $z^0(i)\sim \normal(0,\sigma^2)$ independent and identically distributed.
This can be regarded as  special case of the compressed sensing
model (\ref{eq:obsdata}), whereby $n=N$ and $A = I$ is the
identity matrix -- i.e. there is no underdetermined system
of equations.  We assume that $x^0$ is sparse. It makes sense to consider
soft thresholding
\[
   \hat{x}^{\tau}(i) = \eta( y(i) ; \tau \sigma ) , \qquad
1 \leq i \leq N,\label{eq:SimpleSoft}
\]
where the soft threshold function (with threshold level $\theta$)
is defined by
\begin{eqnarray}
\eta(x;\theta) = \left\{\begin{array}{ll}
x-\theta & \mbox{ if $\theta<x$,}\\
0     & \mbox{ if $-\theta\le x\le \theta$,}\\
x+\theta  & \mbox{ if $x\le -\theta$.}
\end{array}\right.
\end{eqnarray}
In words, the estimator (\ref{eq:SimpleSoft})
`shrinks' the observations $y$ towards the origin by a multiple
$\tau$ of the
noise level $\sigma$.

In place of studying $\xo$ which are $k$-sparse, \cite{DJHS92,DJ94a}
consider random variables $X$ which obey $\prob\{ X \neq 0 \} \leq \eps$,
where $\eps = k/n$.   So let $\cF_{\eps}$ denote the set of
probability measures
placing all but $\eps$ of their mass at the origin:
\[
   \cF_{\eps}  = \{ \nu\,:\, \nu \mbox{ is probability measure with }
\nu(\{0\}) \ge 1-\eps \}.
\]
We  define the soft thresholding mean square error by
\begin{eqnarray}
\stMSE(\sigma^2;  \nu, \tau)
& \equiv & \E\Big\{\big[\eta\big(X+ {\sigma} \cdot Z;
\tau\sigma\big)-X\big]^2\Big\}\, .\label{eq:stMSEDef}
\end{eqnarray}
Here expectation is with respect to independent random variables
$Z \sim  \normal(0,1)$ and $X\sim \nu$.

It is important to allow general $\sigma$ in calculations below.
However, note to the scale invariance
\begin{eqnarray}
\stMSE(\sigma^2;  \nu, \tau)
= \sigma^2\stMSE(1;  \nu^{1/\sigma}, \tau)
 \, ,
\end{eqnarray}
where $\nu^a$ is the probability distribution obtained by rescaling $\nu$:
$\nu^a(S) = \nu(\{x:\, a\, x\in S\})$.
It follows that all calculations can be made in the $\sigma=1$ setting and
results rescaled to obtain final answers.
Below, when we deal with $\sigma=1$,
we will suppress the $\sigma$ argument, and simply write
$\stMSE(\nu, \tau) \equiv \stMSE(1;  \nu, \tau)$

\begin{figure}
\center{\includegraphics[width=10.cm]{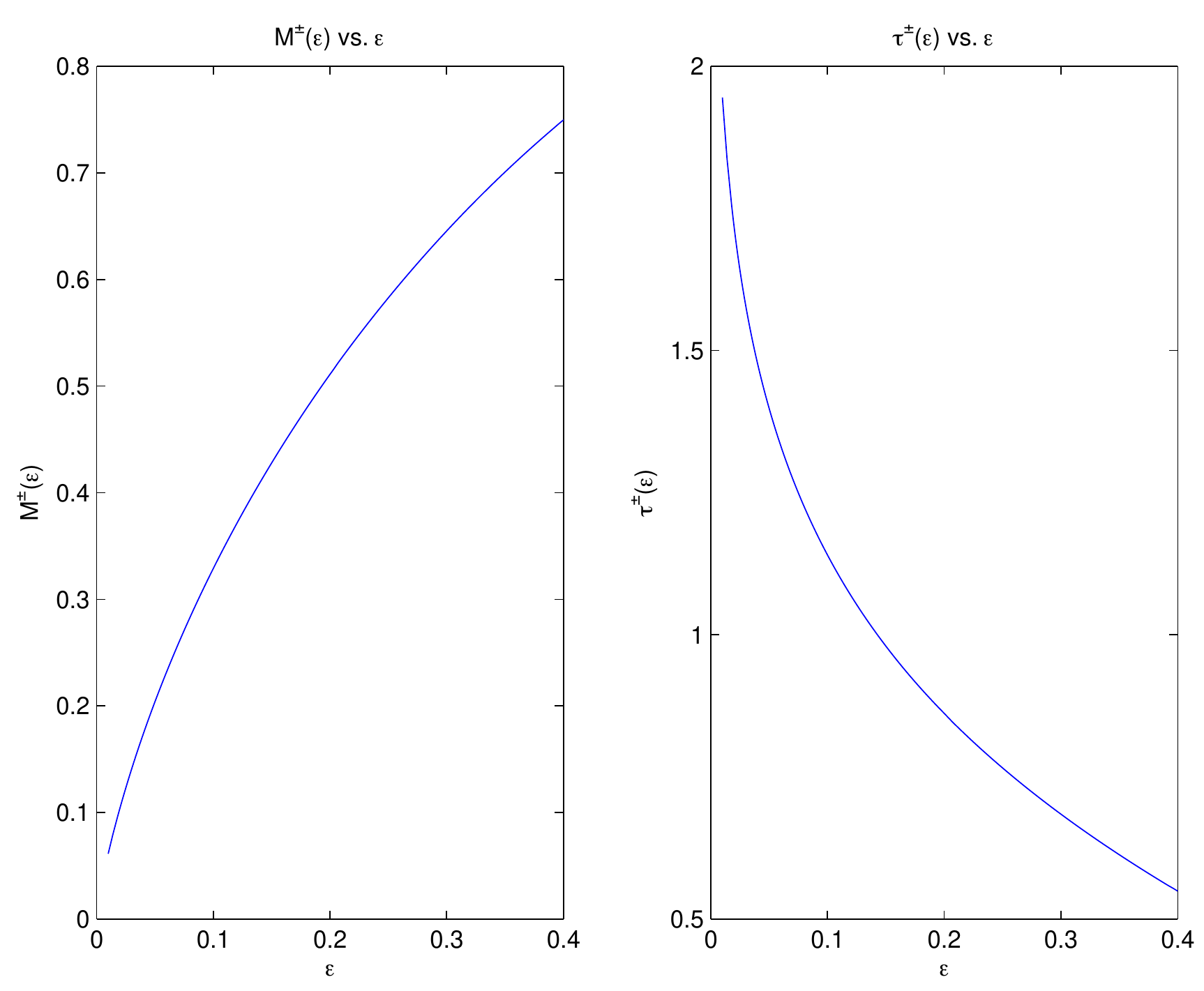}}
\caption{Left:  $M^\hash(\eps)$ as a function of $\eps$;
Right: $\tau^\hash(\eps)$ as a function of $\eps$.}
\label{fig:scalarMMX}
\end{figure}

The {\it minimax threshold MSE} was defined  in \cite{DJHS92,DJ94a}   by
\begin{eqnarray}
    M^\hash(\eps) = \inf_{\tau > 0 } \sup_{\nu \in \cF_\eps}
\stMSE(\nu,\tau)\, .\label{eq:FirstMMAX}
\end{eqnarray}
(The superscript $\pm$ reminds us that, when the estimand $X$ is nonzero,
it may take either sign. In Section \ref{sec-Positivity}, the superscript $+$ will
be used to cover the case where $X \geq 0$).
We will denote by $\tau^{\hash}(\eps)$ the threshold level achieving the
infimum.
Figure \ref{fig:scalarMMX} depicts the behavior of
$M^\hash$  and $\tau^\hash$ as a function of $\eps$.
$M^\hash(\eps)$  was studied in \cite{DJ94a} where one can
find a considerable amount of information
about the behavior of the optimal threshold $\tau^\hash$ and the least
favorable distribution $\nu^\hash_{\eps}$.
In particular, the optimal threshold behaves as
\[
   \tau^\hash(\eps) \sim \sqrt{ 2 \log( \eps^{-1} )}\, , \qquad \mbox{ as }\;\;
\eps \goto 0,
\]
and is explicitly computable at finite $\eps$.

A peculiar aspect of the results in \cite{DJ94a} requires us to generalize their results somewhat.
For a given, fixed $\tau > 0$, the worst case MSE obeys
\begin{eqnarray}
      \sup_{\nu \in \cF_\eps}  \stMSE(\nu,\tau) =
       \eps\, (1 + \tau^2) + (1-\eps)[2(1+\tau^2)\, \Phi(-\tau)-2\tau\,\phi(\tau)]   \, ,\label{eq:SupExpression}
\end{eqnarray}
with $\phi(z)= \exp(-z^2/2)/\sqrt{2\pi}$ the standard normal density
and $\Phi(z) = \int_{-\infty}^z\phi(x)\, \de x$ the Gaussian distribution.
 This supremum is ``achieved'' only
by a three-point mixture on the {\it extended} real line $\reals \cup \{ -\infty , \infty \}$:
\[
    \nu^*_\eps  = (1-\eps) \delta_0 + \frac{\eps}{2} \delta_{\infty} + \frac{\eps}{2} \delta_{-\infty}.
\]
We will need approximations which place no mass at $\infty$.
We say  distribution $\nu_{\eps,\alpha}$ is
{\it $\alpha$-least-favorable} for $ \eta(\,\cdot\,;\tau)$
if it is the least-dispersed distribution in $\cF_\eps$ achieving  a fraction
$(1-\alpha)$
of the worst case risk for $\eta(\,\cdot\,;\tau)$, i.e. if both $(i)$
\[
       \stMSE(\nu_{\eps,\alpha},\tau^{\hash}(\eps))  = (1 -\alpha) \cdot
\sup_{\nu \in \cF_\eps}  \stMSE(\nu,\tau^{\hash}(\eps))\, ,
\]
and $(ii)$ $\nu$ has the smallest second moment for which $(i)$ is true.
The  least favorable distribution $\nu_{\eps,\alpha}$ has the form of
a three-point mixture
\[
    \nu_{\eps,\alpha}  = (1-\eps)\, \delta_0 + \frac{\eps}{2} \delta_{\mu^\hash(\eps,\alpha)} + \frac{\eps}{2} \delta_{-\mu^\hash(\eps,\alpha)}\, .
\]
Here $\mu^\hash(\eps,\alpha)$ is an explicitly computable function,  see below, and for $\alpha > 0$ fixed we have
\[
   \mu^{\hash}(\eps,\alpha) \sim \sqrt{ 2 \log( \eps^{-1} )}\, ,  \qquad \mbox{ as }\;\;
\eps \goto 0\, .
\]
Note in particular the relatively weak role played by $\alpha$. This shows that although the precise least-favorable
situation places mass at infinity, in fact, an approximately least-favorable situation is already achieved much closer to
the origin.

\begin{figure}
\center{\includegraphics[width=10.cm]{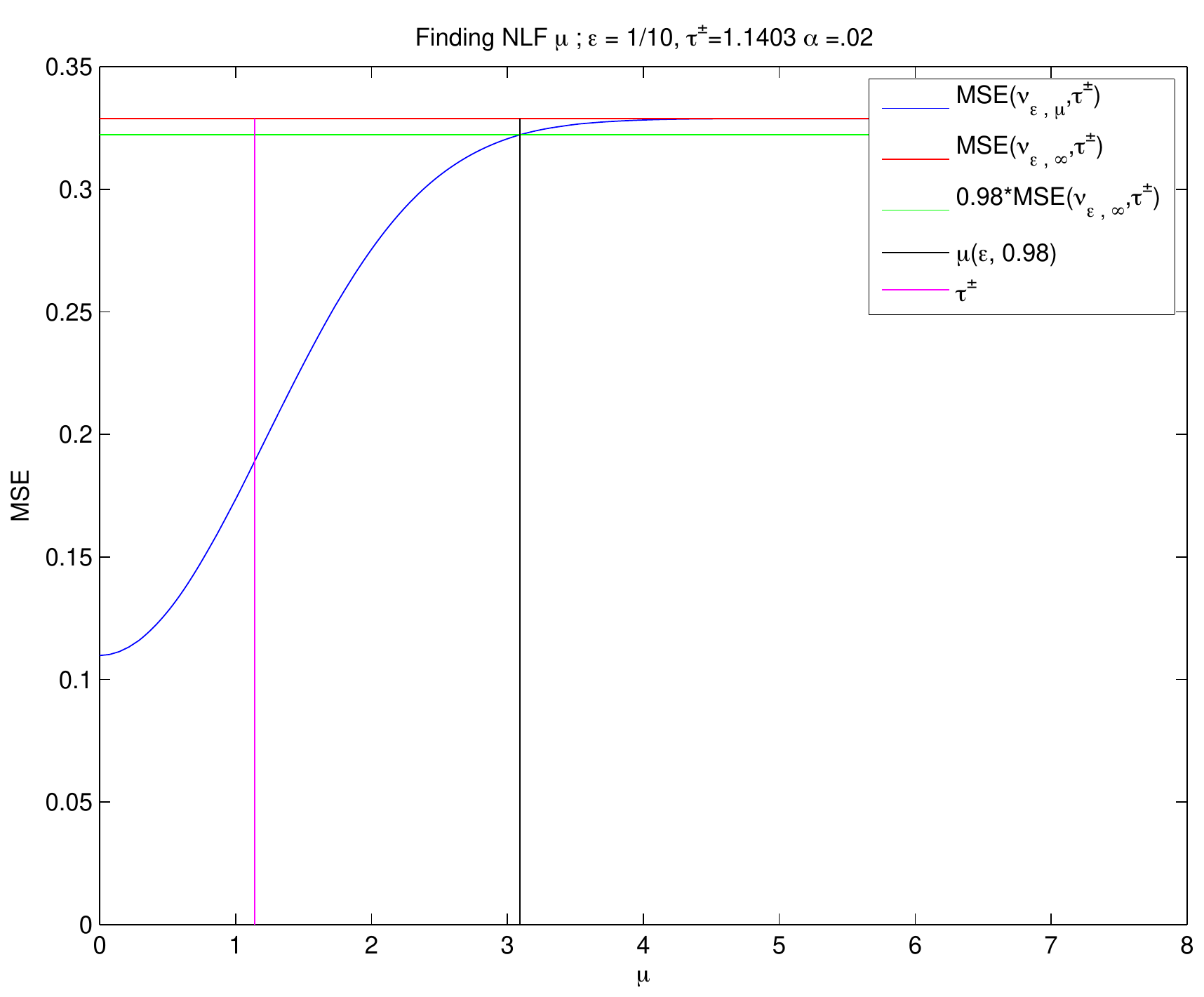}}
\caption{Illustration of $\alpha$-least-favorable $\nu$.
For $\eps=1/10$, we consider
soft thresholding with  the minimax parameter $\tau^{\hash}(\eps)$.
We identify the smallest $\mu$ such that the measure
$\nu_{\eps,\mu} = (1-\eps) \delta_0 + \frac{\eps}{2} \delta_{\mu} +
\frac{\eps}{2} \delta_{-\mu}$ has  $\stMSE(\nu_{\eps,\mu},\tau^*) \geq
0.98\, M^{\hash}(0.1)$ (i.e. the MSE is at least $98\,\%$ of the minimax
MSE).}
\label{fig:scalarLF}
\end{figure}

%
%
\section{Main Results}\label{sec:MainResults}
\newcommand{\hx}{\hat{x}}
\newcommand{\cW}{\cal W}

The notation of the last section allows us to state
our main results.

\subsection{Terminology}

\begin{definition} {\bf (Large-System Limit)}.
A sequence of problem size parameters $n,N$ will be said to {\bf grow proportionally}
if both $n,N \goto \infty$ while $n/N \goto \delta \in (0,1)$.

Consider a sequence of random variables $(W_{n,N})$, where $n,N$ grow proportionally.
Suppose  that $W_{n,N}$ converges in probability to a deterministic quantity $W_\infty$,
which may depend on $\delta > 0$.
Then we say that  $W_{n,N}$ has {\bf large-system limit} $W_\infty$, denoted
\[
        W_\infty = \lslim(W_{n,N}) .
\]
\end{definition}

\begin{definition}
{\bf (Large-System Framework)}.
We denote by $\LSF(\delta,\rho,\sigma,\nu)$
a sequence of  problem instances $(y,A,\xo)_{n,N}$ as per
Eq.~(\ref{eq:obsdata})
indexed by problem sizes $n,N$ growing proportionally:
$n/N \goto \delta$. In each instance, the entries of the $n \times  N$ matrix
$A$ are Gaussian iid $\normal(0, 1/n)$,
the entries of $z^0$ are Gaussian iid $\normal(0,\sigma^2)$ and
the entries of $x^0$ are iid $\nu$.
\end{definition}

For the sake of concreteness we focus here on problem sequences
whereby the matrix $A$ has iid Gaussian entries. An obvious generalization
of this setting would be to assume that the entries are iid
with mean $0$ and variance $1/n$. We expect our result to hold for a broad
set of distributions in this class.

In order to match the $k$-sparsity condition underlying (\ref{eq:obsdata})
we consider the standard framework only for $\nu \in \cF_{\delta\rho}$.
\begin{definition} {\bf (Observable).}
Let $\hx$ denote the output of a reconstruction algorithm on problem instance $(y,A,\xo)$.
An {\em observable}  $J$ is a function $ J(y,A,x^0,\hx)$ of the tuple
$(y,A,x^0,\hx)$.
\end{definition}

In an abuse of notation, the realized values $J_{n,N} = J(y,A,x^0,\hx)$ in this framework
will also be called observables.
An example is the observed per-coordinate
MSE:
\[
\MSE  \equiv\frac{1}{N}\, \| \hx -  x^0 \|_2^2\, .
\]
The $\MSE$ depends explicitly on $x^0$
and implicitly on $y$ and $A$ (through the reconstruction algorithm).
Unless specified, we shall assume that the reconstruction algorithm
solves the LASSO problem (\ref{BPDN}), and hence $\hxl = \hx$.
Further in the following we will drop the dependence
of the observable on the arguments $y,A,x^0,\hx$, and
the problem dimensions $n,N$, when clear from context.
\begin{definition} {\bf  (Formalism).}
A \emph{formalism} is a procedure that assigns
a {\em purported} large-system limit $\Formal(J)$  to  an observable $J$
in the $\LSF(\delta,\rho,\sigma,\nu)$.   This limit in general depends  on $\delta$, $\rho$,
$\sigma^2$, and $\nu \in \cF_{\delta\rho}$: $\Formal(J) = \Formal(J; \delta,\rho,\sigma,\nu)$.
\end{definition}
Thus, in sections below we will consider $J = \MSE(y,A,\xo, \hxl)$
and describe a specific formalism
yielding  $\Formal(\MSE)$, the formal MSE (also denoted by $\FMSE$).
Our formalism has the following character when applied to MSE:
for each $\sigma^2$, $\delta$, and
probability measure $\nu$ on $\reals$,
it calculates
a purported limit $\FMSE(\delta,\nu,\sigma)$.
For a problem instance with large $n,N$
realized from the standard framework $\LSF(\delta,\rho,\sigma,\nu)$,
we claim the MSE will be approximately $\FMSE(\delta,\nu,\sigma)$ .
In fact we will show how to calculate formal limits  for several observables.
For clarity, we always attach the modifier \emph{formal}
to any result of our formalism: e.g., \emph{formal $\MSE$},
\emph{formal False Alarm Rate}, \emph{formally optimal threshold parameter},
and so on.

\begin{definition} {\bf (Validation).}
A formalism is \emph{theoretically} validated by proving that,
in the standard asymptotic framework,
we have
\[
     \lslim(J_{n,N}) = \Formal(J)
\]
for a class $\cJ$ of observables to which the formalism applies,
and for a range of $\LSF(\delta,\rho,\sigma^2,\nu)$.

A formalism is \emph{empirically} validated by showing that,
for problem instances $(y,A,\xo)$ realized from
$\LSF( \delta,\rho,\sigma,\nu)$
with large $N$ we have
\[
 J_{n,N} \approx \Formal(J; \delta,\rho,\sigma,\nu),
\]
for a collection of observables $J \in \cJ$ and a range of asymptotic
framework parameters $( \delta,\rho,\sigma,\nu)$; here the
approximation $\approx$ should be evaluated by usual
standards of empirical science.
\end{definition}
Obviously, theoretical validation is stronger than empirical
validation, but careful empirical validation is still validation.
We do not attempt here to theoretically validate
this formalism in any generality;
see \cite{BayatiMontanari} results in this direction.
Instead we view the formalism as calculating {\it predictions}
of empirical results.
We have compared these predictions with empirical results and found
a persuasive level of agreement.
For example, our formalism has been used to predict the
MSE of reconstructions by (\ref{BPDN}), and actual empirical
results match the predictions, i.e.:
\[
 \frac{1}{N}  \| \hxl - x^0 \|_2^2\approx \FMSE(\delta,\rho,\nu,\sigma) .
\]
%
%
\subsection{Results of the Formalism}

The behavior of formal mean square error changes dramatically at the
following phase boundary.
\begin{definition}[Phase Boundary]
For each $\delta \in [0,1]$, let $\rhoMSE(\delta)$
be the value of $\rho$ solving
\begin{eqnarray}
M^\hash(\rho\delta) = \delta\, .
\end{eqnarray}
\end{definition}
It is well known that $M^\hash(\eps)$ is monotone increasing
and concave in $\eps$,
with $M^{\hash}(0)=0$ and $M^{\hash}(1)=1$. As a consequence,
$\rhoMSE$ is also a monotone increasing function of $\delta$,
in fact $\rhoMSE(\delta) \goto 0$ as $\delta \goto 0$
and $\rhoMSE(\delta) \goto 1$ as $\delta \goto 1$.
An explicit expression for the curve $(\delta,\rhoMSE(\delta))$
is provided in Appendix \ref{app:Formulae}.
%
%

\begin{proposition} \label{prop:main} {\bf Results of Formalism.}
The formalism developed below yields the following conclusions.

\begin{description}
\item{1.a} In the region $\rho < \rhoMSE(\delta)$, the minimax formal noise
sensitivity obeys the formula
\[
   M^*(\delta, \rho) \equiv \frac{M^\hash(\rho\delta)}{1 - M^\hash(\rho\delta)/\delta} .
\]
In particular, $M^*$ is finite throughout this region.

\item{1.b} With $\sigma^2$ the noise level in (\ref{eq:obsdata}),
define the {\em formal noise-plus interference} level\phantom{,}
$\FNPI = \FNPI(\tau;\delta,\rho,\sigma,\nu) $
\[
\FNPI=   \sigma^2 + \FMSE/\delta,
\]
and its minimax value $\NPI^*(\delta,\rho; \sigma) \equiv \sigma^2 \cdot ( 1 + M^*(\delta,\rho)/\delta)$.
For $\alpha > 0$,  define
\[
 \mu^*(\delta,\rho; \alpha) \equiv \mu^{\hash}(\delta\rho,\alpha) \cdot  \sqrt{\NPI^*(\delta,\rho)}
\]
In $\LSF(\delta,\rho,\sigma,\nu)$
let $\nu \in \cF_{\delta \rho}$ place fraction $1-\delta\rho$ of its mass at zero and the remaining mass equally on
$\pm \mu^*(\delta,\rho;\alpha)$. This $\nu$ is
$\tilde{\alpha}$-least-favorable:
the formal noise sensitivity of $\hxl$ equals
$(1-\tilde{\alpha}) M^*(\delta,\rho)$, with
$(1-\tilde{\alpha}) =(1-\alpha)(1-M^{\hash}(\delta\rho))/(1-(1-\alpha)
M^{\hash}(\delta\rho))$.

\item{1.c} The {\em formally maximin penalty} parameter obeys
\[
\lambda^*(\nu; \delta,\rho, \sigma) \equiv   \tau^\hash(\delta\rho) \cdot  \sqrt{\FNPI(\tau^\hash;\delta,\rho,\sigma,\nu)}
\cdot ( 1 - \EqDR(\nu; \tau^\hash(\delta\rho))/\delta)\, ,
\]
where $\EqDR(\,\cdots\,)$ is the asymptotic \emph{detection rate},
i.e. the asymptotic fraction of coordinates that are
estimated to be nonzero.
(An explicit expression for this quantity is given in Section
\ref{eq:EqDR_formula}.)

In particular with this $\nu$-adaptive choice of penalty parameter, the formal
$\MSE$ of $\hxl$  does not exceed $M^* \cdot \sigma^2$.

\item{2} In the region $\rho > \rhoMSE(\delta)$, the formal noise
sensitivity  is infinite. Throughout this phase, for each fixed number
$M < \infty$, there exists $\alpha > 0$  such that
the probability distribution $\nu \in \cF_{\delta\rho}$ placing
its nonzeros at $\pm  \mu^*(\delta,\rho,\alpha)$,
yields formal $\MSE$ larger than $M$.
\end{description}
\end{proposition}

We explain the formalism and derive these results in Section \ref{sec:Formalism} below.

\subsection{Interpretation of the Predictions}

Figure \ref{fig:Figure001} displays the noise sensitivity; above the
phase transition boundary $\rho = \rhoMSE(\delta)$, it is infinite.
The different contour lines show positions in the $\delta,\rho$
plane where a given noise sensitivity is achieved. As one might expect,
the sensitivity blows up rather dramatically as we approach the
phase boundary.

\begin{figure}
\center{\includegraphics[width=10.cm]{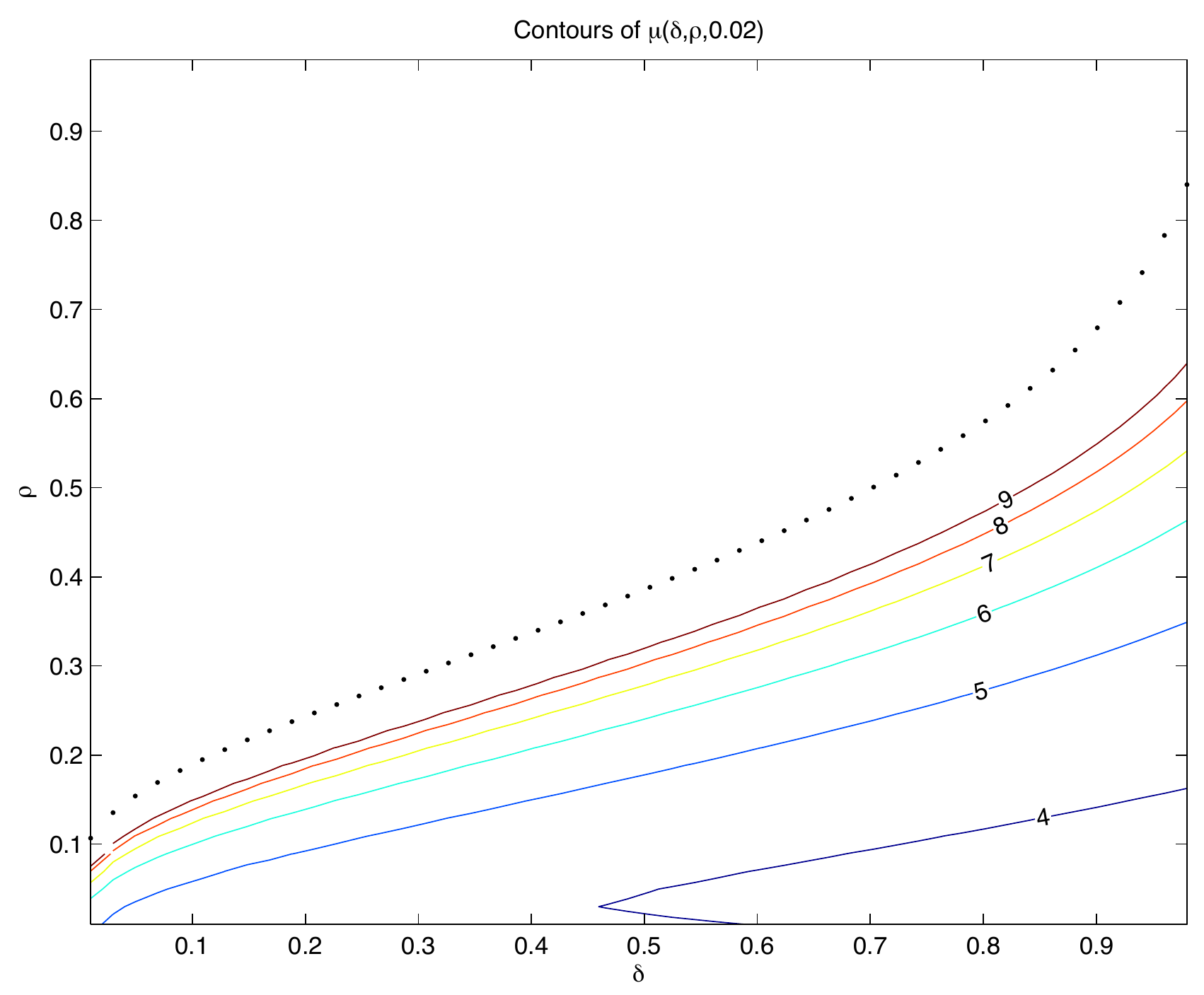}}
\caption{\small Contour lines of the near-least-favorable signal amplitude
$\mu^*(\delta, \rho,\alpha)$ in the $(\rho,\delta)$ plane. The
dotted line corresponds to the phase transition $(\delta,\rhoMSE(\delta))$,
while the colored solid lines portray level sets of $\mu^*(\delta, \rho,\alpha)$. The 3-point mixture distribution
$(1-\eps) \delta_0 + \frac{\eps}{2} \delta_{\mu} + \frac{\eps}{2} \delta_{-\mu}$, $(\eps = \delta\rho)$ will cause
98\% of the worst-case MSE. When a $k$-sparse vector is drawn
from this distribution, its nonzeros are all at $\pm \mu$.}
\label{fig:AltFigure002}
\end{figure}

Figure \ref{fig:AltFigure002} displays the least-favorable coefficient
amplitude $\mu^*(\delta,\rho,\alpha=0.02)$. Notice that
$\mu^*(\delta,\rho,\alpha)$ diverges as the phase boundary is approached.
Indeed beyond the phase boundary arbitrarily large MSE can be
produced by choosing $\mu$ large enough.

\begin{figure}
\center{\includegraphics[width=10.cm]{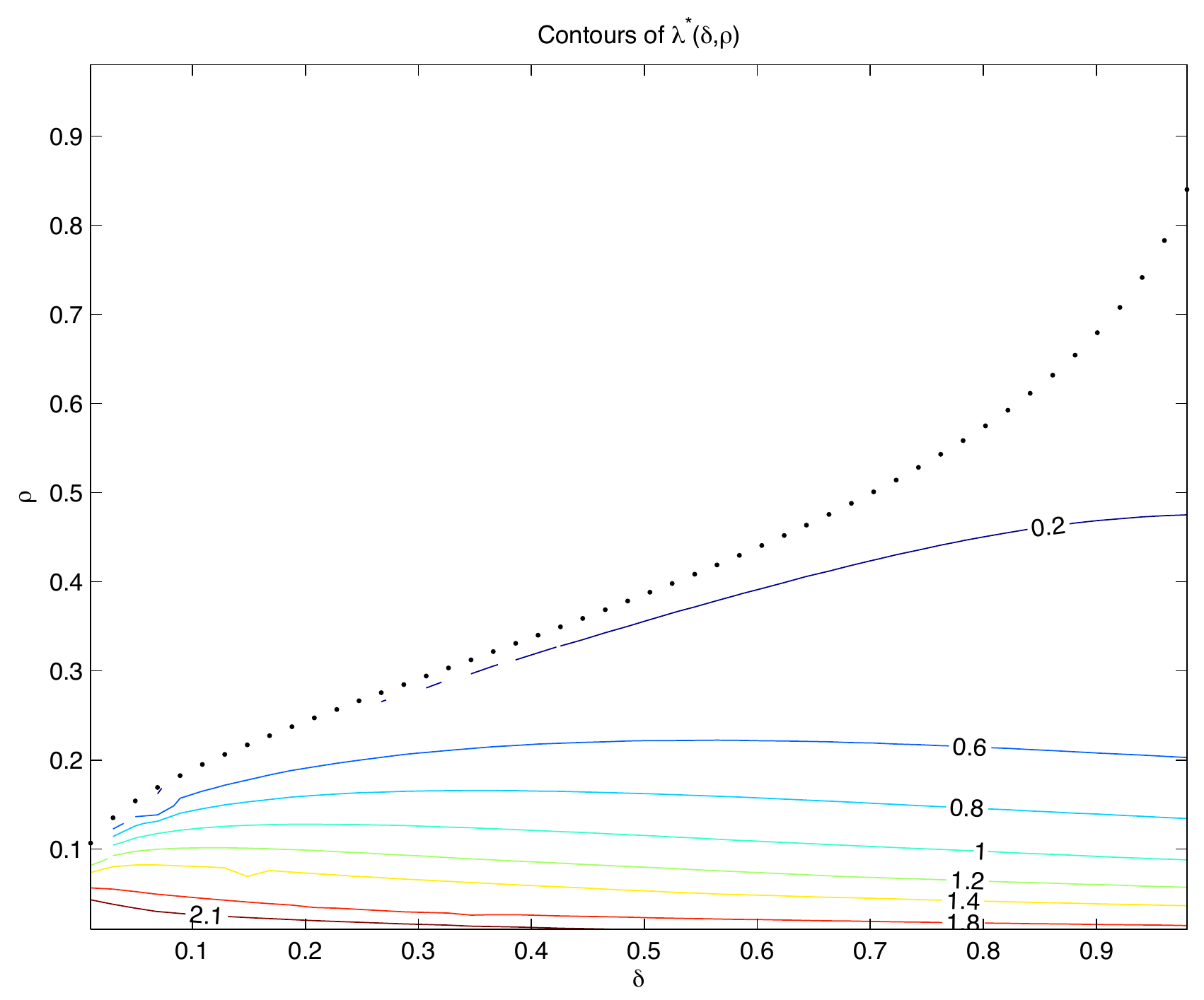}}
\caption{\small Contour lines of the maximin penalization parameter:
$\lambda^*(\delta, \rho)$ in the $(\rho,\delta)$ plane. The
dotted line corresponds to the phase transition $(\delta,\rhoMSE(\delta))$,
while thin lines are contours for $\lambda^*(\delta, \rho,\alpha)$. Close to phase
transition, the maximin value approaches $0$.}
\label{fig:ContoursLambdaStar}
\end{figure}

Figure \ref{fig:ContoursLambdaStar} displays the value of
the optimal penalization parameter
amplitude $\lambda^* = \lambda^*(\nu_{\delta,\rho}^*; \delta,\rho,\sigma=1)$. Note that the
parameter tends to zero as we approach phase transition.

For these figures, the region above phase transition
is not decorated, because the values there are infinite or not defined.

\subsection{Comparison to other phase transitions}

In view of the importance of the phase boundary for Proposition
\ref{prop:main}, we note the following:
\begin{finding}  {\bf Phase Boundary Equivalence.} \label{prop:bdryequiv}
The phase boundary $\rhoMSE$ is identical to the
phase boundary $\rho_{\ell_1}$ below which $\ell_1$ minimization
and $\ell_0$ minimization are equivalent.
\end{finding}

In words, throughout the phase where $\ell_1$ minimization is equivalent to
$\ell_0$ minimization, the solution to (\ref{BPDN}) has bounded formal MSE.
When we are outside that phase, the solution has unbounded formal MSE.
The verification of Finding \ref{prop:bdryequiv} follows in two steps.
First, the formulas for the phase boundary discussed in this paper
are identical to the phase boundary formulas given in \cite{AMPSupplement};
Second, in \cite{AMPSupplement} it was shown that these
formulas agree numerically with the formulas known for $\rho_{\ell_1}$.

\subsection{Validating the Predictions}
Proposition \ref{prop:main} makes predictions for the behavior of  solutions to (\ref{BPDN}).
It will be validated empirically, by showing that such solutions behave as predicted.

In particular, simulation evidence will be presented to show that in the
phase where noise sensitivity is {\it finite}:
\begin{enumerate}
\item Running (\ref{BPDN}) for data $(y,A)$ generated from vectors
$x_0$ with coordinates with distribution $\nu$ which is nearly least-favorable
results in an empirical MSE approximately equal to
$M^*(\delta,\rho) \cdot \sigma^2$.
\item Running (\ref{BPDN})  for data $(y,A)$ generated from vectors
$x_0$ with coordinates with distribution $\nu$ which is
far from least-favorable
results in empirical MSE noticeably smaller than $M^*(\delta,\rho) \cdot \sigma^2$.
\item Running (\ref{BPDN})  with a suboptimal penalty parameter $\lambda$
results in empirical MSE noticeably greater than $ M^*(\delta,\rho)\cdot  \sigma^2$.
\end{enumerate}
Second,  in the phase where formal MSE is {\it infinite}:
\begin{enumerate}
\item[4.] Running (\ref{BPDN}) on vectors $x_0$ generated by formally least-favorable
results in an empirical MSE which is very large.
\end{enumerate}
Evidence for all these claims will be given below.
%
%
\section{The formalism}\label{sec:Formalism}

\subsection{The AMPT Algorithm}

We now consider a reconstruction approach seemingly very different
from ($\BPDN$). This algorithm, called
\emph{first-order approximate message passing}
(AMP) algorithm proceeds iteratively, starting at $\bhat^0 = 0$ and
producing the estimate $\bhatee$ of $\targ$ at iteration $t$
according to the iteration:
\begin{eqnarray}
z^t & = & y -A\bhatee +z^{t-1} \frac{\df_{t}}{n}
\label{eq:FOAMP2} \\
\bhatpo & = & \eta(A^*z^t+\bhatee;\theta_t)\, ,
\label{eq:FOAMP1}
\end{eqnarray}
Here  $\bhatee\in\reals^p$ is the current estimate of $\targ$, and $\df_t = \| \bhatee \|_0$
is the number of nonzeros in the current estimate. Again
$\eta(\,\cdot\,; \,\cdot\,)$  is the {\it soft threshold} nonlinearity
with threshold parameter $\theta_t$
\begin{equation} \label{ThreshChoice}
   \theta_t = \tau \cdot \sigma_t ;
\end{equation}
$\tau$ is a tuning constant, fixed throughout iterations
and $\sigma_t$ is an empirical measure of the scale of the residuals.
Finally $z^{t}\in\reals^n$ is the current {\it working residual}.
Compare with the usual residual defined by $r^t = y - A \bhatee$
via the identity $z^t = r^t  + z^{t-1} \frac{\df_{t}}{n}$.
The extra term in AMP plays a subtle but crucial role.
\footnote{A similar-looking algorithm was introduced by the authors
in \cite{DMM}, with identical steps (\ref{eq:FOAMP1})-(\ref{eq:FOAMP2});
it differed only in the choice of threshold;
instead of a tuning parameter  $\tau$
like in (\ref{ThreshChoice}) -- one that can be set freely -- a
 fixed choice $\tau(\delta)$ was made for each specific $\delta$.
Here we call that algorithm AMPM \--  $M$  for {\it minimax},
as explained in  \cite{AMPSupplement}.  In contrast,
the current algorithm is tunable, allowing choice of $\tau$,
we label it AMPT$(\tau)$, $T$ for tunable.}
%
%
\subsection{Formal MSE, and its evolution}

Let $\npi(m;\sigma,\delta) \equiv \sigma^2 + m/\delta$. We define the
\emph{MSE map} $\Psi$ through
\begin{eqnarray}
 \Psi(m,\delta,\sigma,\tau,\nu)
\equiv  \stMSE(\npi(m,\sigma,\delta); \nu,\tau)\, ,\label{eq:PsiDef}
\end{eqnarray}
where the function $\stMSE(\,\cdot\,; \nu,\tau)$ is the soft thresholding
mean square error already introduced in Eq.~(\ref{eq:stMSEDef}).
It describes the MSE of soft thresholding in a problem where the noise level
is $\sqrt{\npi}$. A heuristic explanation of the meaning and origin of
$\npi$ will be given below.

\begin{definition}
{\bf State Evolution}.
{\sl The {\em state} is a 5-tuple $(m; \delta, \sigma,\tau,\nu) $.
{\em State evolution} is  the evolution of the state by the rule
\begin{eqnarray*}
(m_t; \delta, \sigma, \tau,\nu) &\mapsto& (\Psi(m_t); \delta,  \sigma,
\tau,\nu),\\
 t  &\mapsto&  t+1 .
\end{eqnarray*}
As the parameters  $(\delta,\sigma,\tau, \nu)$ remain
fixed during evolution, we usually omit mention of them and think of
state evolution simply as the  iterated application of $\Psi$:
\begin{eqnarray*}
m_t&\mapsto& m_{t+1} \equiv \Psi(m_t),\\
 t  &\mapsto&  t+1 .
\end{eqnarray*}
}
\end{definition}

\begin{figure}
\center{\includegraphics[width=10.cm]{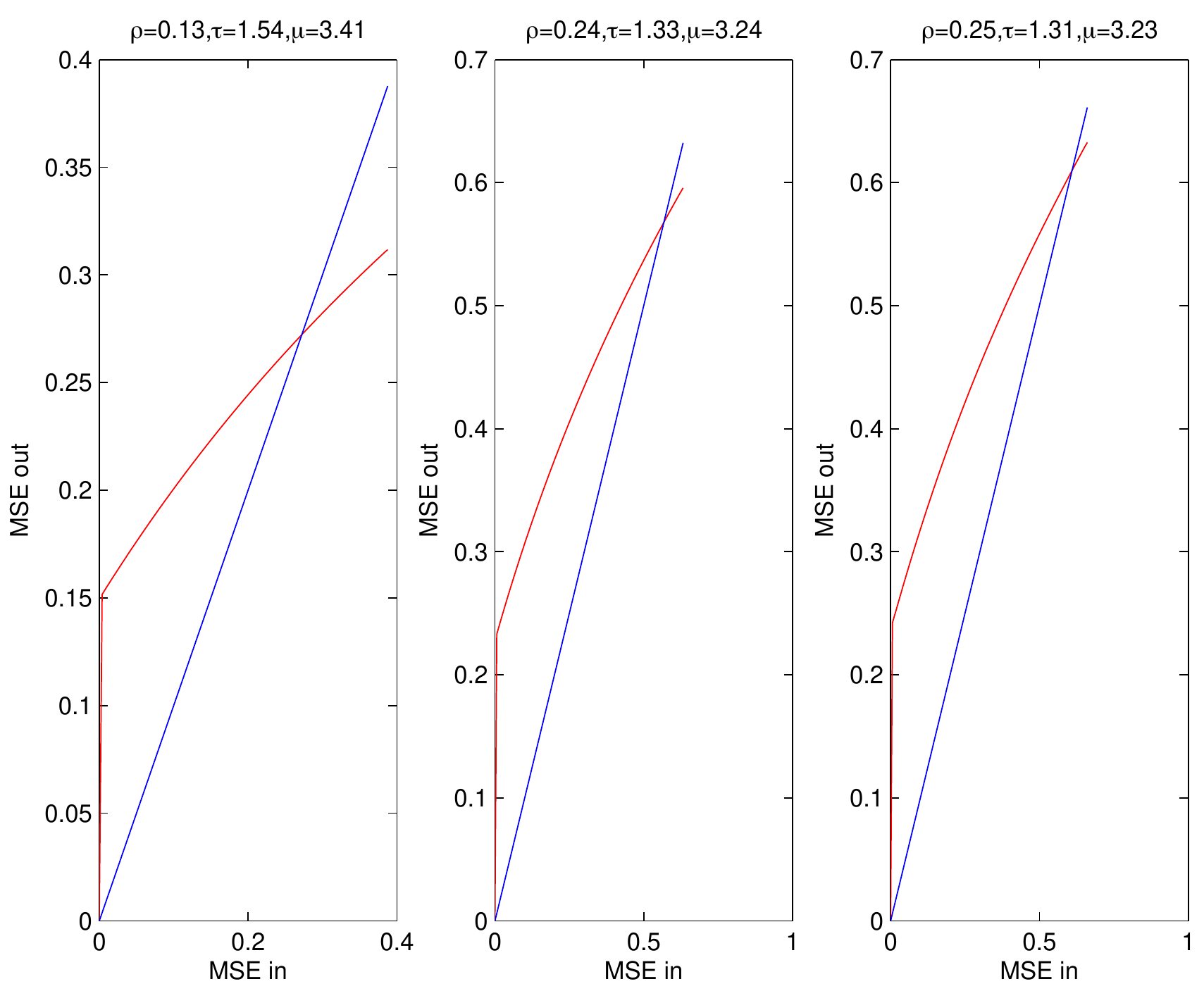}}
\caption{MSE Map $\Psi$ in three cases, and associated fixed points.
Left:  $\delta=0.25$, $\rho = \rhoMSE/2$, $\sigma=1$, $\nu = \nu^*(\delta,\rho,\alpha)$
Center:  $\delta=0.25$, $\rho = \rhoMSE \times 0.95$, $\sigma=1$, $\nu = \nu^*(\delta,\rho,\alpha)$
Right: $\delta=0.25$, $\rho = \rhoMSE$, $\sigma=1$, $\nu = \nu^*(\delta,\rho,\alpha)$ }
\label{fig:MSEMapsA}
\end{figure}

\begin{definition}
{\bf Stable Fixed Point}.
The Highest Fixed Point of the continuous function $\Psi$ is
\[
    \HFP(\Psi) = \sup\{ m : \Psi(m) \geq m \}.
\]
The stability coefficient of the continuously differentiable function
$\Psi$ is
\[
    \SC(\Psi) = \left.\frac{\de}{\de m} \Psi(m) \right|_{m = \HFP(\Psi)} .
\]
We say that $\HFP(\Psi)$ is a stable fixed point if $0 \leq \SC(\Psi) < 1$.
\end{definition}
To illustrate this, Figure \ref{fig:MSEMapsA} shows the MSE
map and fixed points in three cases.

In what follows we denote by $\mu_2(\nu) = \int x^2 \de\nu$ the
second-moment of the distribution $\nu$.
\begin{lemma}\label{lemma:Convergence}
Let $\Psi(\,\cdot\,) = \Psi(\,\cdot\,,\delta,\sigma,\tau,\nu)$,
and assume either $\sigma^2>0$ or $\mu_2(\nu)>0$.
Then the sequence of iterates $m_t$ defined by  $m_{t+1}= \Psi(m_{t})$
starting from $m_0 = \mu_2(\nu)$
converges monotonically to $\HFP(\Psi)$:
\[
           m_t \goto \HFP(\Psi) , \qquad t \goto \infty.
\]
Further, if $\sigma>0$ then $\HFP(\Psi)\in (0,\infty)$ is the unique fixed point.

Suppose further that the stability coefficient satisfies
$0 < \SC(\Psi) < 1$. Then there exists a constant ${\cal A}(\nu,\Psi)$
such that
\[
  \big|m_t - \HFP(\Psi) \big| \leq {\cal A}(\nu,\Psi)\, \SC(\Psi)^t\, .
\]
Finally, if $\mu_2(\nu) \ge \HFP(\Psi)$ then the sequence $\{m_t\}$
is monotonically decreasing to $\mu_2(\nu)$ with
\[
 (m_t - \HFP(\Psi)) \leq \SC(\Psi)^t  \cdot (\mu_2(\nu) - \HFP(\Psi)) .
\]
\end{lemma}
In short, barring the trivial case $x^0=0$, $z^0=0$
(no signal, no noise),
state evolution converges to the highest fixed point.
If the stability coefficient is smaller than $1$, convergence is
exponentially fast.

\begin{proof}[Proof (Lemma \ref{lemma:Convergence})]
 This Lemma is an immediate consequence
of the fact that $m\mapsto \Psi(m)$ is a concave non-decreasing function,
with $\Psi(0)>0$ as long as $\sigma>0$ and $\Psi(0) = 0$ for
$\sigma = 0$.

Indeed in \cite{AMPSupplement} the authors showed that at noise level
$\sigma=0$, the MSE map  $m\to\Psi(m ;\delta,\sigma,\nu,\tau)$
is concave as a function of $m$.
We have the identity
\[
    \Psi( m ;\delta,\sigma,\nu,\tau) =  \Psi( m  + \sigma^2 \cdot \delta ;\delta,\sigma=0,\nu,\tau),
\]
relating the noise-level $0$ MSE map to the noise-level $\sigma$ MSE map.
From this it follows that $\Psi$ is concave for $\sigma>0$
as well. Also,  \cite{AMPSupplement} shows that
$\Psi( m=0  ;\delta,\sigma=0,\nu,\tau)=0$ and
$\frac{\de\Psi}{\de m}( m=0  ;\delta,\sigma=0,\nu,\tau)>0$,
whence $\Psi( m=0  ;\delta,\sigma,\nu,\tau)>0$ for any positive noise level
$\sigma$.
\end{proof}

In the same paper \cite{AMPSupplement},
the authors derived the least-favorable
stability coefficient in the noiseless case $\sigma=0$:
\[
 \SC^*(\delta,\rho, \sigma=0) =  \sup_{\nu \in \cF_{\delta\rho}}
\SC(\Psi(\,\cdot\, ;\delta,\sigma=0,\nu,\tau))\, .
\]
They showed that, for $M^{\hash}(\delta,\rho)<\delta$
the only fixed point is at $m=0$ and has stability coefficient
\[
     \SC^*(\delta,\rho, \sigma=0)  =  M^\hash(\delta\rho)/\delta\, .
\]
Hence, it follows that $\SC^*(\delta,\rho, \sigma=0)  < 1$ throughout
the region $\rho < \rhoMSE(\delta)$.

Define
\[
   \SC^*(\delta,\rho) =  \sup_{\sigma > 0} \sup_{\nu \in \cF_{\delta\rho}}
\SC(\Psi(\, \cdot\, ;\delta,\sigma,\nu,\tau)).
\]
Concavity of the noise level 0 MSE map implies
\[
 \SC^*(\delta,\rho)  = \SC^*(\delta,\rho, \sigma=0)) .
\]
We therefore conclude that throughout the region $\rho < \rhoMSE(\delta)$
For this reason, that region can also be called the {\it  stability phase},
not only the stability coefficient is smaller than $1$,
$\SC(\Psi)<1$, but that it can be bounded away from $1$ uniformly
in the signal distribution $\nu$.
\begin{lemma}
Throughout the region $\rho < \rhoMSE(\delta)$, $0 < \delta < 1$,
for every $\nu \in \cF_{\delta\rho}$,
we have
$\SC(\Psi) \leq \SC^*(\delta,\rho) < 1$.
\end{lemma}

\begin{figure}
\center{\includegraphics[width=10.cm]{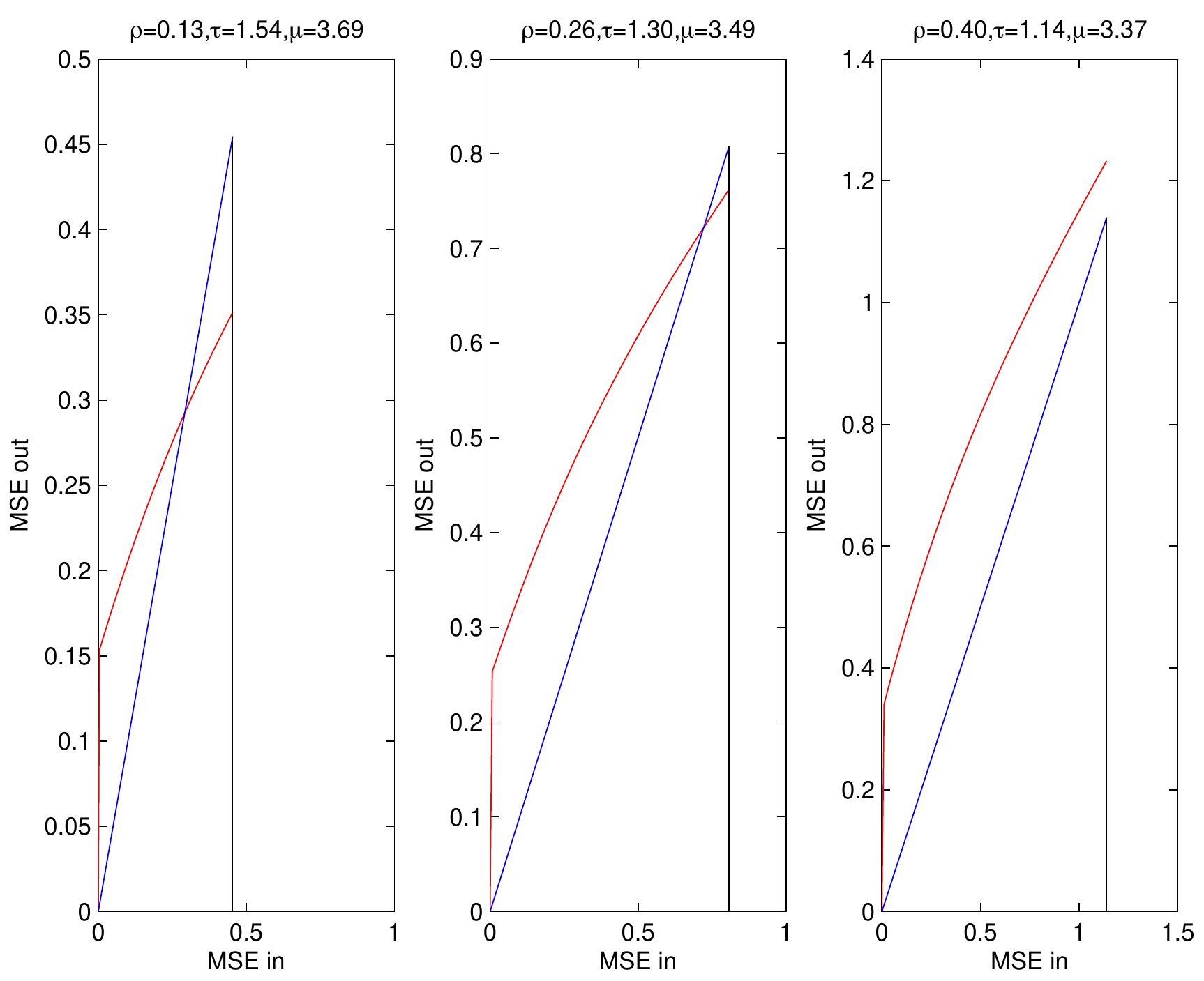}}
\caption{Crossing the phase transition: effects on MSE Map $\Psi$, and associated state evolution.
Left:  $\delta=0.25$, $\rho = \rhoMSE/2$, $\sigma=1$, $\nu = \nu(\delta,\rho,0.01)$
Middle: $\delta=0.25$, $\rho = 0.9 \cdot \rhoMSE$, $\sigma=1$, $\nu = \nu(\delta,\rho,0.01)$
Right: $\delta=0.25$, $\rho = 1.5 \cdot \rhoMSE$, $\sigma=1$, $\nu = \nu(\delta,\rho,0.01)$.
In each case $\tau = \tau^{\hash}(\delta\rho)$.}
 \label{fig:MSEMapsB}
\end{figure}

Outside the stability region, for each large $m$,
we can find measures $\nu$ obeying the sparsity constraint $\nu \in \cF_{\delta\rho}$ for which
state evolution converges to a fixed point suffering equilibrium $\MSE > m$. The construction
in section \ref{sec-AbovePT} shows that
$\HFP(\Psi) > \mu_2(\nu) > m $.
Figure \ref{fig:MSEMapsB} shows the MSE
map and the state evolution in three cases which may be compared to \ref{fig:MSEMapsA}.
In the first case, $\rho$ is well below $\rhoMSE$ and the fixed point is well below $\mu_2(\nu)$.
In the  second case, $\rho$ is slightly below $\rhoMSE$ and the fixed point is
close to $\mu_2(\nu)$. In the third case, $\rho$ is above $\rhoMSE$ and the fixed point,
lies above $\mu_2(\nu)$.

$\mu_2(\nu)$ is the MSE one suffers by `doing nothing':
setting threshold $\lambda = \infty$ and taking $\hx = 0$.  When $\HFP(\Psi) > \mu_2(\nu)$,
one iteration of thresholding makes things {\it worse}, not better.  In words, the phase boundary is exactly
 the place below which we are sure that, if $\mu_2(\nu)$ is large, a single iteration of thresholding gives an estimate $\hx^1$ that is better
than the starting point $\hx^0$.  Above the phase boundary, even a single iteration of thresholding
may be a catastrophically bad thing to do.

\begin{definition}
{\bf (Equilibrium States and State-Conditional Expectations)}

Consider a real-valued function $\zeta: \reals^3 \mapsto \reals$,
its \emph{expectation in state $S =  (m; \delta,\sigma,\nu)$} is
\[
      \cE(\zeta | S) = \E\,\big\{ \zeta(X,Z,
\eta(X+\sqrt{\npi}\, Z;\tau\sqrt{\npi})) \big\}\, ,
\]
where $\npi = \npi(m;\sigma,\delta)$ and $X \sim \nu$,  $Z \sim \normal(0,1)$
are independent random variables.

Suppose we are given $(\delta,\sigma,\nu,\tau)$,
and a fixed point $m^*$, $m^* = \HFP(\Psi)$
with $\Psi = \Psi(\, \cdot\, ; \delta,\sigma,\nu,\tau)$.
The tuple $S^*= (m^*; \delta,\sigma,\nu)$ is called the {\em equilibrium state}
of state evolution.
The {\em expectation in the equilibrium state}
is $ \cE(\zeta | S^*)$.
\end{definition}

\begin{definition}  {\bf (State Evolution Formalism for AMPT) }.
Run the AMPT algorithm and assume that the sequence of
estimates $(\hx^t,z^t)$ converges to the fixed point
$(\hx^{\infty}, z^{\infty})$.  To each function
$\zeta: \reals^3 \mapsto \reals$ associate the
observable
\[
J^\zeta ( y,A,\xo,\hx) = \frac{1}{N}
\sum_{i=1}^N \, \zeta\big( \xo(i) , A^Tz(i)+\hx(i)-\xo(i) , \hx(i)\big)  \, .
\]
Let $S^*$ denote the equilibrium state reached by state evolution
in a given situation $(\delta,\sigma,\nu,\tau)$.
The state evolution formalism assigns the purported limit value
\[
      \Formal(J^\zeta) =  \cE(\zeta | S^*).
\]
\end{definition}

Validity of the state evolution formalism for AMPT
entails that, for a sequence of problem instances $(y,A,\xo)$
drawn from $\LSF(\delta,\rho,\sigma,\nu)$,
 the large-system limit for observable $ J^\zeta_{n,N} $
is simply the expectation in the equilibrium state:
\[
   \lslim J^\zeta_{n,N} =  \cE(\zeta | S^*).
\]

The class $\cJ$ of observables representable by the form
$J^\zeta$ is quite rich, by choosing $\zeta(u,v,w)$ appropriately. Table \ref{table:observables} gives
examples of well-known observables and
the $\zeta$ which will generate them.
\begin{table}
\begin{center}
\begin{tabular}{|l|l|l|}
\hline
Name   & Abbrev. & $\zeta=\zeta(u,v,w)$ \\
\hline
Mean Square Error  &MSE   & $\zeta = (u-w)^2$     \\
False Alarm Rate     & FAR         & $\zeta = 1_{\{  w \neq 0 \& u =  0\}}/(1-\rho\delta)$  \\
Detection Rate          & DR                &
$\zeta = 1_{\{  w \neq 0\}}$ \\
Missed Detection Rate  &MDR & $\zeta = 1_{\{ w =  0 \& u \neq 0  \}}/(\rho\delta)$ \\
False Detection Rate  & FDeR  & $\zeta = 1_{\{ w \ne  0 \& u =  0  \}}/(\rho\delta)$ \\
\hline
\end{tabular}
\caption{Some observables and their names.}
\label{table:observables}
\end{center}
\end{table}
Formal values for other interesting observables can in principle be obtained
by combining such simple ones.  For example, the False Discovery rate FDR
is the ratio FDeR$/$DR and so the ratio of two elementary observables
of the kind for which the formalism is defined.
We assign it the purported limit value
\[
     \Formal(\FDR) =  \frac{\Formal(\FDeR)}{\Formal(\DR)}\, .
\]
Below we list a certain number of observables
for which the formalism was checked empirically and that play an important
role in characterizing the fixed point estimates.

\vspace{0.5cm}

\begin{centering}
{\bf Calculation of Formal Operating Characteristics of $\AMPT(\tau)$ by State Evolution}
\end{centering}

\bitem
 \item[] Given $\delta, \sigma, \nu,\tau$, identify the fixed point
$\HFP(\Psi(\,\cdot\, ;\delta,\sigma,\nu,\tau )$.
Calculate the following quantities
 \bitem
 \item Equilibrium MSE
 \[
        \EqMSE = m_\infty = \HFP(\Psi(\,\cdot\, ;\nu,\tau); \delta,\sigma).
 \]
 \item Equilibrium Noise Plus Interference Level
 \[
      \npi_\infty =  \frac{1}{\delta}m_\infty + \sigma^2
 \]
 \item Equilibrium Threshold (absolute units)
 \[
      \theta_\infty = \tau \cdot \sqrt{\npi_\infty}.
 \]
  \item Equilibrium Mean Squared Residual. Let $Y_\infty = X +
\sqrt{\npi_\infty}\, Z$
for $X \sim \nu$ and $Z \sim \normal(0, 1)$ are independent. Then
 \[
     \EqMSR =
\E  \big\{\left[ Y_\infty - \eta(Y_\infty ; \theta_\infty) \right]^2 \big\}\, .
 \]
  \item Equilibrium Mean Absolute Estimate
 \[
        \EqMAE =  \E \{| \eta(Y_\infty ; \theta_\infty) |\}\, .
 \]
  \item Equilibrium Detection Rate
\begin{eqnarray}
        \EqDR =  \prob \{  \eta(Y_\infty  ; \theta_\infty) \neq 0 \}\, .
\label{eq:EqDR_formula}
\end{eqnarray}
  \item Equilibrium Penalized MSR
 \[
        \EqPMSR = \EqMSR/2 +
\theta_\infty \cdot ( 1 - \EqDR/\delta) \cdot \EqMAE .
 \]
\eitem

\eitem


\vspace{.1in}



\noindent
%
%
\subsection{AMPT - LASSO Calibration}

Of course at this point the reader is entitled to feel that the
introduction of AMPT is a massive digression.  The relevance of AMPT
is indicated by the following conclusion from \cite{TheoPredLassoOpChar}:
\begin{finding} \label{find:crude}
In the large system limit, the operating characteristics
of $\AMPT(\tau)$ are
equivalent to those of \LASSO$(\lambda)$
under an appropriate calibration $\tau \leftrightarrow \lambda$.
\end{finding}

By {\it calibration}, we mean a rescaling that maps results
on one problem into results on the other problem.  The notion
is explained at greater length in \cite{TheoPredLassoOpChar}.
The correct mapping can be guessed from the following remarks:
\bitem
\item[]  $\LASSO(\lambda)$: no {\it residual}
exceeds $\lambda$:
$\|A^T(y - A \hxl)\|_{\infty} \leq \lambda$.
Further
\begin{eqnarray*}
\hxl_i > 0 &\Leftrightarrow& (A^T(y - A \hxl))_i = \lambda\, ,\\
\hxl_i = 0 &\Leftrightarrow& |(A^T(y - A \hxl))_i| < \lambda\, ,\\
\hxl_i < 0 &\Leftrightarrow& (A^T(y - A \hxl))_i = -\lambda\, .
\end{eqnarray*}
\item  $\AMPT(\tau)$:  At a fixed point
$\bhat^{\infty}$, $z^{\infty}$, no {\it working residual}
exceeds the equilibrium threshold $\theta_\infty$:
$\| A^Tz^{\infty}\|_{\infty}\le\theta_{\infty}$. Further
\begin{eqnarray*}
\bhat^{\infty}_i > 0 &\Leftrightarrow& (A^T z^{\infty})_i = \theta_{\infty}\, ,\\
\bhat^{\infty}_i = 0 &\Leftrightarrow& |(A^T z^{\infty})_i| < \theta_{\infty}\, ,\\
\bhat^{\infty}_i < 0 &\Leftrightarrow& (A^T z^{\infty})_i = -\theta_{\infty}\, .
\end{eqnarray*}
\eitem
Define $df = \#\{i : \bhat^{\infty}_i \neq 0\}$.
Further notice that at the AMPT fixed point $(1-\df/n)z^{\infty}=
y-A^T\bhat^{\infty}$.
We can summarize these remarks in the following statement
\begin{lemma}
Solutions $\hxl$ of  $\LASSO(\lambda)$ (i.e. optima of the problem
(\ref{BPDN})) are in correspondence with fixed points
$(\bhat^{\infty},z^{\infty})$ of the
$\AMPT(\tau)$ under the bijection $\bhat^{\infty}=\hxl$,
$z^{\infty}=  (y-A^T\hxl)/(1-\df/n)$, provided the threshold parameters
are in the following relation
\begin{eqnarray}
   \lambda = \theta_\infty \cdot ( 1 - \df/n)\, .\label{eq:Calibration}
\end{eqnarray}
\end{lemma}
In other words, if we have a fixed point of $\AMPT(\tau)$
we can choose $\lambda$ in such a way that this is also an optimum of
$\LASSO(\lambda)$. Viceversa, any optimum of $\LASSO(\lambda)$
can be realized as a fixed point of $\AMPT(\tau)$: notice in fact that
the relation (\ref{eq:Calibration}) is invertible
whenever $\df<n$.

This simple rule gives a calibration
relationship between $\tau$ and $\lambda$, i.e.
a one-one correspondence between $\tau$ and $\lambda$
that renders the two apparently different
reconstruction procedures equivalent, provided the iteration
$\AMPT(\tau)$ converges rapidly to its fixed point.
Our empirical results confirm that this is indeed
what happens for typical
large system frameworks $\LSF(\delta,\rho,\sigma,\nu)$.

The next lemma characterizes the equilibrium calibration relation
between AMP and LASSO.
\begin{lemma}
Let $\EqDR(\tau) =
\EqDR(\tau;\delta,\rho,\nu,\sigma)$ denote the equilibrium detection rate
obtained from state evolution when the tuning parameter of AMPT is $\tau$.
Define $ \tau^0(\delta,\rho,\nu,\sigma) > 0$, so that
$\EqDR(\tau) \leq \delta$ when $\tau > \tau^0$. For each $\lambda \geq 0$,
there is a unique value $\tau(\lambda) \in  [\tau_0,\infty)$
such that
\[
        \lambda = \theta_\infty(\tau) \cdot ( 1 - \EqDR(\tau) /\delta).
\]
\end{lemma}

We can restate Finding \ref{find:crude} in the following more convenient
form.
\begin{finding} \label{find:fine}
For each $\lambda \in [0,\infty)$
we find that $\AMPT(\tau(\lambda))$ and  $\LASSO(\lambda)$
have statistically equivalent observables. In particular the $\MSE$,
$\MAE$, $\MSR$, $\DR$,
have the same distributions.
\end{finding}
%
%
\subsection{Derivation of Proposition \ref{prop:main}}
\label{sec:MainPropo}

Consider the following Minimax Problem for $\AMPT(\tau)$.
With $\FMSE(\tau;\delta,\rho,\sigma,\nu)$ denoting the equilibrium
formal MSE for
$AMPT(\tau)$ for the framework $\LSF(\delta,\rho,\sigma,\nu)$, fix $\sigma=1$
and \emph{define}
\begin{equation} \label{eq:AltMstarDef}
    M^\flat(\delta,\rho) = \inf_\tau \sup_{\nu \in \cF_{\delta \rho}}
\FMSE(\tau;\delta,\rho,\sigma=1,\nu).
\end{equation}

We will first show that this definition
obeys the formula just like the one in Proposition \ref{prop:main},
given for $M^*$.
Later we show that $M^\flat = M^*$.
\begin{proposition} \label{prop:MMxFormula} For $M^\flat$ defined by (\ref{eq:AltMstarDef}),
\begin{equation} \label{MinimaxFormula}
    M^\flat(\delta,\rho) =  \frac{ M^\hash(\delta\rho) }{1- M^\hash(\delta\rho)/\delta}
\end{equation}
The AMPT threshold rule
\begin{equation} \label{MinimaxThresh}
     \tau^*(\delta,\rho) = \tau^\hash(\delta\rho), \quad 0 < \rho < \rhoMSE(\delta)\, ,
\end{equation}
minimaxes the formal MSE:
\begin{equation} \label{eq:ThreshMinMax}
 \sup_{\nu \in \cF_{\delta \rho}} \FMSE(\tau^*;\delta,\rho,1,\nu)   =
\inf_\tau \sup_{\nu \in \cF_{\delta \rho}} \FMSE(\tau;\delta,\rho,1,\nu)   = M^\flat(\delta,\rho).
\end{equation}
\end{proposition}

\begin{figure}
\center{\includegraphics[width=10.cm]{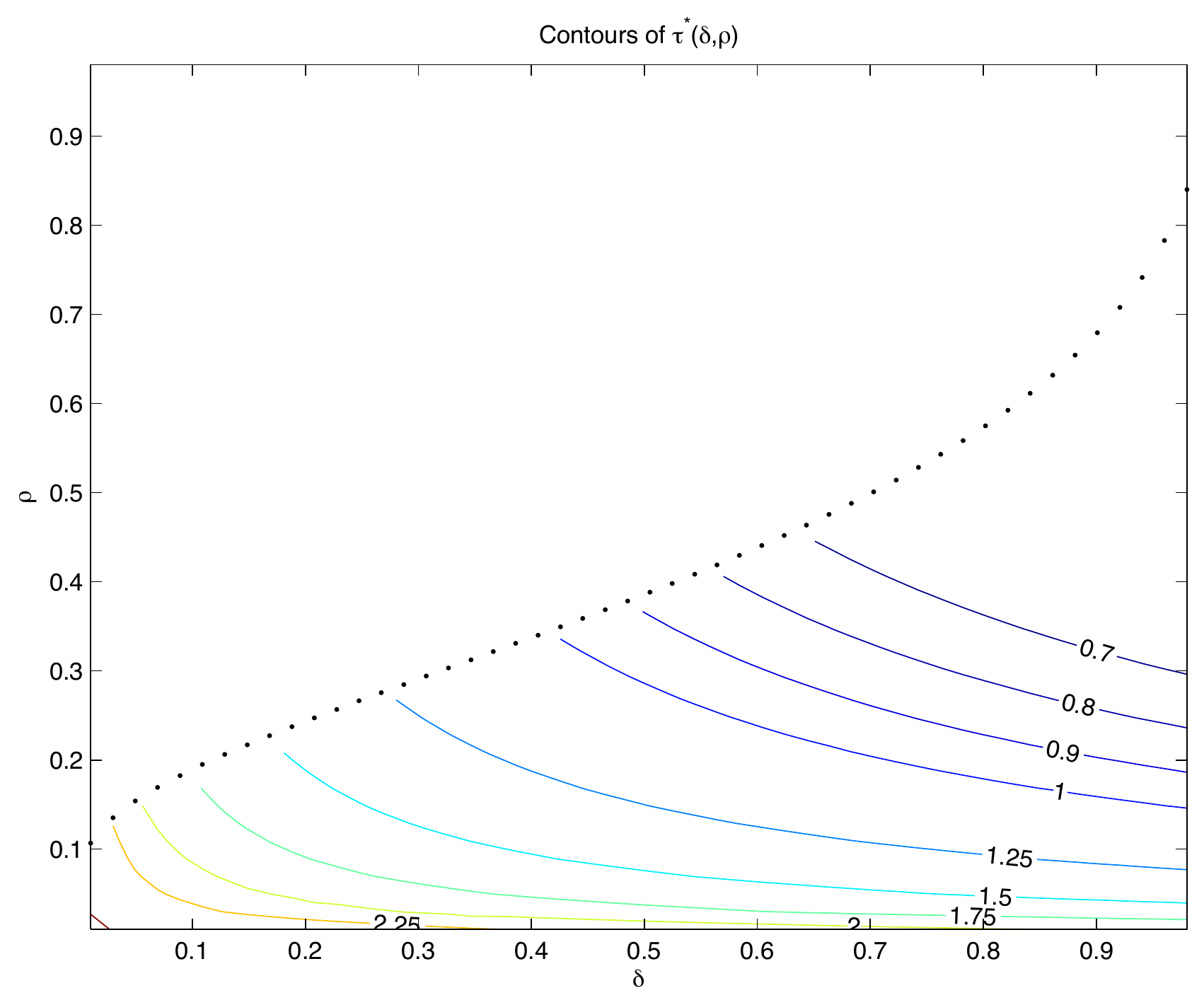}}
\caption{\small Contour lines of $\tau^*(\delta,\rho)$ in the $(\rho,\delta)$ plane. The
dotted line corresponds to the phase transition $(\delta,\rhoMSE(\delta))$,
while thin lines are contours for $\tau^*(\delta, \rho)$}
\label{fig:ContoursTauStar}
\end{figure}

Figure \ref{fig:ContoursTauStar} depicts the behavior of $\tau^*$ in the $(\delta,\rho)$ plane.

\begin{proof}[Proposition \ref{prop:MMxFormula}]
Consider $\nu \in \cF_{\delta\rho}$ and $\sigma^2=1$  and set
$\tau^*(\delta,\rho) = \tau^\hash(\delta\rho)$ as in the statement.
Let for short
$\Psi(m;\nu) =  \Psi(m,\delta,\sigma=1,\tau^*,\nu)=
\stMSE(\npi(m,1,\delta) ; \nu,\tau^*)$, cf. Eq.~(\ref{eq:PsiDef}).
Then $m^* = \HFP(\Psi)$ obeys, by definition of fixed point,
\[
     m^* = \Psi(m^*; \nu)  \, .
\]
 We can use the scale invariance
$\stMSE(\sigma^2 ; \nu,\tau^*) = \stMSE(1 ; \tilde{\nu},\tau^*)$,
where $\tilde{\nu}$ is a rescaled probability measure, $\tilde{\nu}\{ x \cdot \sigma \in B \} = \nu \{ x \in B \}$.
For $\nu \in \cF_{\delta\rho}$, we have $\tilde{\nu} \in \cF_{\delta\rho}$
as well and we therefore obtain
\[
   m^* = \stMSE(\npi(m^*,1,\delta) ; \nu,\tau^*)
=  \stMSE(1 ; \tilde{\nu},\tau^*)\cdot\npi(m^*,1,\delta)
\leq M^\hash(\delta\rho) \cdot  \npi(m^*;1, \delta)\, \, ,
\]
where we used the fact that $\tau^*(\delta,\rho) = \tau^\hash(\delta\rho)$.
Hence
\[
      \frac{m^*}{\npi(m^*;1,\delta) } \leq M^\hash(\delta\rho)\,  .
\]
The function $m \mapsto \frac{m}{\npi(m;\delta,1)}$ is one-to-one
strictly increasing on the
interval $[0,\delta)$.
Thus, provided that $1-  M^\hash(\delta\rho) /\delta > 0$, i.e. $\rho < \rhoMSE$, we have
\[
     m^* \leq \frac{ M^\hash(\delta\rho)}{ 1-  M^\hash(\delta\rho) /\delta}.
\]
As this inequality applies to {\it any} $\HFP$ produced by our formalism, in particular the
largest one consistent with $\nu \in \cF_{\delta\rho}$, we have
\[
   \sup_{\nu \in \cF_{\delta \rho}} \FMSE(\tau^*;\delta,\rho,1,\nu)   \leq \frac{ M^\hash(\delta\rho)}{ 1-  M^\hash(\delta\rho) /\delta}.
\]

We now develop the reverse inequality.
To do so, we make a specific choice $\onu$ of $\nu$.
Fix $\alpha > 0$ small.
Now for $\eps = \delta\rho$, define
$\xi = \mu^{\hash}(\eps,\alpha) \cdot \sqrt{\NPI^*}$,
where $\NPI^* = 1 + M^\flat/\delta$ (with $M^\flat =
M^\hash(\delta\rho)/(1-M^\hash(\delta\rho)/\delta)$
as in the thesis).
Let $\overline{\nu} = (1-\eps) \delta_0 + (\eps/2)\, \delta_{-\xi} +
(\eps/2) \delta_\xi$. Denote by $m^* = m^*(\onu)$ the highest fixed point
corresponding to the signal distribution $\onu$.
 Using once again scale invariance, we have
\begin{eqnarray}
  m^*= \stMSE(\npi(m^*,1,\delta) ; \onu,\tau^*)
=  \stMSE(1 ; \tilde{\nu},\tau^*)\cdot\npi(m^*,1,\delta) \, ,
\label{eq:RelationAgain}
\end{eqnarray}
where $\tilde{\nu}$ is again a  rescaled probability measure,
this time with $\tilde{\nu}\{ x \cdot \sqrt{\npi(m^*,1,\delta)} \in B \} = \overline{\nu} \{ x \in B \}$.
Now since $m^* \leq M^\flat$, we have
$\npi(m^*,1,\delta)  \leq \NPI^*$, and hence
\[
     \frac{\xi}{\sqrt{\npi(m^*,1,\delta)}} =\mu^{\hash}(\eps,\alpha)  \cdot  \sqrt{\frac{\NPI^*}{\npi(m^*,1,\delta)}}  > \mu^{\hash}(\eps,\alpha)\,  .
\]
Note that $\stMSE(m;(1-\eps) \delta_0 + (\eps/2) \delta_{-x} +
(\eps/2) \delta_x ,\tau)$ is monotone increasing in $|x|$.
Recall that  $\nu_{\eps,\alpha} = (1-\eps) \delta_0 + (\eps/2)
\delta_{-\mu^{\hash}(\eps,\alpha)} + (\eps/2) \delta_{\mu^{\hash}
(\eps,\alpha)}$
is $\alpha$-least favorable for the minimax problem (\ref{eq:FirstMMAX}).
Consequently,
\[
    \stMSE(1 ; \tilde{\nu},\tau^*) \geq
\stMSE(1 ; \nu_{\delta\rho,\alpha},\tau^*)  = (1-\alpha) \cdot
M^\hash(\delta,\rho)\, .
\]
Using the scale-invariance relation, Eq.~(\ref{eq:RelationAgain}),
we conclude that
\[
      \frac{m^*}{\npi(m^*;\delta,1) } \geq (1-\alpha) \cdot M^\hash(\delta\rho)\, .
\]
Again, in the region $\rho < \rhoMSE(\delta)$,  the function $m \mapsto \frac{m}{\npi(m;\delta,1)}$ is one-to-one and monotone and therefore
so
\[
    \FMSE(\tau^*;\delta,\rho,1,\onu)    \geq \frac{ (1-\alpha) \cdot M^\hash(\delta\rho)}{ 1-  (1-\alpha) \cdot M^\hash(\delta\rho) /\delta}.
\]
As $\alpha > 0$ is arbitrary, we conclude
\[
    \sup_{\nu \in \cF_{\delta \rho}}
\FMSE(\tau^*;\delta,\rho,1,\nu) \geq \frac{ M^\hash(\delta\rho)}{ 1-  M^\hash(\delta\rho) /\delta}.
\]
\end{proof}

We now explain how this result about AMPT leads to our claim for the
behavior of the LASSO estimator $\hxl$.
By a scale invariance the quantity (\ref{eq:DefMStar}) can be rewritten
as a fixed-scale $\sigma=1$ property:
\[
M^*(\delta,\rho) = \sup_{\nu \in \cF_{\delta \rho}}\inf_\lambda
\FMSE(\nu,\lambda | \LASSO)\, ,
\]
where we introduced explicit reference to the algorithm used,
and dropped the irrelevant arguments.
We will analogously write $\FMSE(\nu,\tau | \AMPT)$
for the AMPT$(\tau)$ MSE.
\begin{proposition}
 Assume the validity of our calibration relation i.e. the equivalence of formal operating characteristics of $\AMPT(\tau)$
and $\LASSO(\lambda(\tau))$.
Then
\[
    M^*(\delta,\rho) =  M^\flat(\delta,\rho).
\]
Also, for $\lambda^*$ as defined in Proposition \ref{prop:main},
\[
M^*(\delta,\rho) =  \sup_{\nu \in \cF_{\delta \rho}} \FMSE(\nu,\lambda^*(\nu;\delta,\rho,\sigma) | \LASSO )   .
\]
\end{proposition}

In words, $\lambda^*$ is the maximin penalization and
 the maximin MSE of \LASSO is precisely given by the formula (\ref{MinimaxFormula}).

\begin{proof} Taking the validity of our calibration relationship $\tau \leftrightarrow \lambda(\tau)$
as given, we must have
\[
      \FMSE(\nu,\lambda(\tau) | \LASSO ) = \FMSE(\nu,\tau | \AMPT)\, .
 \]
Our definition of $\lambda^*$ in Proposition \ref{prop:main} is simply the calibration relation
applied to the minimax AMPT threshold  $\tau^*$, i.e. $\lambda^* = \lambda(\tau^*)$.
Hence assuming the validity of our calibration relation, we have:
\begin{eqnarray}
          \sup_{\nu \in \cF_{\delta \rho}} \FMSE(\nu,\lambda^*(\nu;\delta,\rho,\sigma) | \LASSO)   &= &
          \sup_{\nu \in \cF_{\delta \rho}}  \FMSE(\nu,\lambda(\tau^*) | \LASSO )  \nonumber \\
           &=&     \sup_{\nu \in \cF_{\delta \rho}}  \FMSE(\nu,\tau^* | \AMPT)  \nonumber  \\
      &=&  \sup_{\nu \in \cF_{\delta \rho}}\inf_\tau  \FMSE(\nu,\tau | \AMPT)  \label{EqMinMax} \\
      &=&  \sup_{\nu \in \cF_{\delta \rho}}\inf_\tau  \FMSE(\nu,\lambda(\tau) | \LASSO)  \nonumber \\
      &=&  \sup_{\nu \in \cF_{\delta \rho}}\inf_\lambda \FMSE(\nu,\lambda | \LASSO)  . \nonumber
 \end{eqnarray}
Display (\ref{EqMinMax}) shows that all these equalities are equal to $M^\flat(\delta,\rho)$.
\end{proof}
The proof of Proposition \ref{prop:main}, points $1a$, $1b$, $1c$
follows immediately from the above.

%
%

\subsection{Formal MSE above Phase Transition}
\label{sec-AbovePT}

We now make an explicit construction showing that noise
sensitivity is unbounded above PT.

We first consider the AMPT algorithm above PT.
Fix $\delta$, $\rho$ with $\rho > \rhoMSE(\delta)$ and
set  $\eps = \delta\rho$.

In this section we focus on 3 point distributions
with mass at $0$ equal to $1-\eps$.
With an abuse of notation we let
$\stMSE(\mu,\tau)$ denote the  MSE of scalar soft thresholding  for
amplitude of the non-zeros equal to $\mu$, and noise variance equal to $1$.
In formulas,
$\stMSE(\mu,\tau)\equiv \stMSE(1;(1-\eps)\delta_0+(\eps/2)\delta_{\mu}+
(\eps/2)\delta_{-\mu},\tau)$, and
\begin{eqnarray*}
\stMSE(\mu,\tau)  = (1-\eps) \E \eta(Z;\tau)^2 + \eps \E \left( \mu - \eta(\mu+Z;\tau) \right)^2\, .
\end{eqnarray*}
Consider values of the AMPT threshold
$\tau$ such that  $\stMSE(0,\tau) < \delta$;
this will be possible for all  $\tau$ sufficiently large.
Pick a number $\gamma \in (0,1)$ obeying
\begin{equation} \label{gammadef}
1 < \gamma < \stMSE(0,\tau)/\delta.
\end{equation}

Let $M^\hash(\eps,\tau)  = \sup_\mu \stMSE(\mu,\tau) $ denote the
worst case risk of $\eta(\,\cdot\, ; \tau)$ over the class $\cF_\eps$.
Let $\mu^{\hash}(\eps,\alpha,\tau)$ denote the $\alpha$-least-favorable $\mu$
for threshold $\tau$:
\[
    \stMSE(\mu^{\hash},\tau) = (1-\alpha) M^\hash(\eps,\tau).
\]
Define $\alpha^* = 1 - \gamma \delta/M^\hash(\eps,\tau)$,
and note that $\alpha^* \in (0,1)$ by earlier assumptions.
Let $\mu^* = \mu^{\hash}(\alpha^*,\tau,\eps)$.
A straightforward calculation along the lines of the previous section yields.
\begin{lemma}
For the measure
$\nu = (1-\eps) \delta_0 + (\eps/2) \delta_{\mu^*} + (\eps/2) \delta_{-\mu^*}$,
the formal $\MSE$ and formal $\NPI$ are given by
\begin{eqnarray*}
\FMSE(\nu,\tau| \AMPT)  &=& \frac{\delta \gamma}{1-\gamma}\, ,\\
\FNPI (\nu,\tau| \AMPT)  &=&  \frac{1}{1-\gamma}\, .
\end{eqnarray*}
\end{lemma}
Assumption (\ref{gammadef}) permits us to choose $\gamma$ very close to 1.
Hence the above formulas show explicitly that MSE is unbounded above phase transition.

What do the formulas say about $\hxl$ above PT?
The $\tau$'s which can be associated to $\lambda$ obey
\[
    0 < \EqDR(\nu,\tau) \leq \delta ,
 \]
where
$\EqDR(\nu,\tau)=\EqDR(\tau;\delta,\rho,\nu,\sigma)$
is the equilibrium detection rate for a signal with distribution
$\nu$.
Equivalently, they are those $\tau$ where the equilibrium discovery
number is $n$ or smaller.
\begin{lemma} \label{lem:LASSOAbovePT}
For each  $\tau > 0$, obeying both
\[
     \stMSE(0,\tau) < \delta \qquad \mbox{ and } \qquad
\EqDR(\nu,\tau)  < \delta,
\]
the parameter $\lambda \geq 0$  defined by  the calibration relation
\[
  \lambda(\tau) = \frac{\tau}{\sqrt{1-\gamma}} \cdot   ( 1  -
\EqDR(\nu,\tau)/\delta),
  \]
has the formal MSE
  \[
        \FMSE(\nu,\tau| \LASSO ) = \frac{\delta \gamma}{1-\gamma}\, .
  \]
\end{lemma}

One can check that, for each $\lambda \geq 0$, for
each phase space point above phase transition,
the above construction allows to construct a
measure $\mu$ with $\eps = \delta\rho$ mass on nonzeros
and with arbitrarily high formal MSE.  This completes the derivation
of part 2 of Proposition \ref{prop:main}.
%
%
\section{Empirical Validation}\label{sec:Empirical}

So far our discussion explains
how state evolution calculations are carried out so others
might reproduce them.  The actual `science contribution' of our paper
comes in showing that these calculations describe the actual
behavior of solutions to (\ref{BPDN}).  We check these calculations
in two ways: first, to show that individual MSE predictions
are accurate, and second, to show that the mathematical structures
(least-favorable, minimax saddlepoint, maximin threshold)
that lead to our predictions are visible in empirical  results.

\subsection{Below phase transition}

Let $\FMSE(\lambda; \delta,\rho,\sigma,\nu)$  denote the formal MSE we assign to $\hxl$
for problem instances from $\LSF(\delta,\rho,\sigma,\nu)$.
Let $\eMSE(\lambda)_{n,N}$
denote the empirical MSE of the LASSO estimator
$\hxl$ in a problem instance drawn from
$\LSF(\delta,\rho,\sigma,\nu)$
at a given problem size $n,N$.  In claiming that the noise sensitivity
of $\hxl$ is bounded above by $M^*(\delta,\rho)$, we are saying that in empirical trials, the  ratio
 $\eMSE/\sigma^2$ will not be larger than $M^*$  with statistical significance.
  We now present empirical evidence for this claim.

\subsubsection{Accuracy of MSE at the LF signal}

We first consider the accuracy of theoretical predictions at the nearly-least-favorable signals
generated by $\nu_{\delta,\rho,\alpha}= (1-\eps) \delta_0 + (\eps/2)
\delta_{-\mu^{*}(\delta,\rho,\alpha)} + (\eps/2) \delta_{\mu^{*}
(\delta,\rho,\alpha)}$
 defined by Part $2.b$ of Proposition \ref{prop:main}.
If the empirical ratio $\eMSE/\sigma^2$ is substantially
above the theoretical bound $M^*(\delta,\rho)$, according to
standards of statistical significance,
we have falsified the proposition.

We consider parameter points $\delta \in \{ 0.10, 0.25,0.50\}$
and $\rho \in \{ \frac{1}{2} \cdot \rhoMSE, \frac{3}{4} \cdot \rhoMSE,
\frac{9}{10} \cdot \rhoMSE, \frac{19}{20} \cdot  \rhoMSE\}$.
The predictions of the SE formalism are detailed in
Table \ref{table:formalPred}.

\begin{table}
\begin{center}
\begin{tabular}{|l|l|l|l|l|l|l|l|l|l|}
\hline
 $\delta$ & $\rho$ & $\eps$ &  $M^\hash(\eps)$ &  $\tau^\hash(\eps)$  & $\mu^\hash(\eps,0.02)$  &$M^*(\delta,\rho)$ & $\mu^*(\delta,\rho,0.02)$ & $\tau^*(\delta,\rho)$& $\lambda^*$\\
\hline
0.10 & 0.09 & 0.01 & 0.06 & 1.96 &  3.74 & 0.14 & 5.79 & 1.96 & 1.28 \\
\hline 0.10 & 0.14 & 0.01 & 0.08 & 1.83 &  3.63 & 0.41 & 8.24 & 1.83 & 0.83 \\
\hline 0.10 & 0.17 & 0.02 & 0.09 & 1.77 &  3.58 & 1.20 & 12.90 & 1.77 & 0.51 \\
\hline 0.10 & 0.18 & 0.02 & 0.10 & 1.75 &  3.57 & 2.53 & 18.28 & 1.75 & 0.41 \\
\hline 0.25 & 0.13 & 0.03 & 0.15 & 1.54 &  3.41 & 0.39 & 5.46 & 1.54 & 0.98 \\
\hline 0.25 & 0.20 & 0.05 & 0.20 & 1.40 &  3.29 & 1.12 & 7.68 & 1.40 & 0.62 \\
\hline 0.25 & 0.24 & 0.06 & 0.23 & 1.33 &  3.24 & 3.28 & 12.22 & 1.33 & 0.39 \\
\hline 0.25 & 0.25 & 0.06 & 0.24 & 1.31 &  3.23 & 6.89 & 17.31 & 1.31 & 0.30 \\
\hline 0.50 & 0.19 & 0.10 & 0.32 & 1.15 &  3.11 & 0.90 & 5.19 & 1.15 & 0.70 \\
\hline 0.50 & 0.29 & 0.14 & 0.42 & 1.00 &  2.99 & 2.55 & 7.35 & 1.00 & 0.42 \\
\hline 0.50 & 0.35 & 0.17 & 0.47 & 0.92 &  2.93 & 7.51 & 11.75 & 0.92 & 0.26 \\
\hline 0.50 & 0.37 & 0.18 & 0.48 & 0.90 &  2.91 & 15.75 & 16.67 & 0.90 & 0.20 \\
\hline
\end{tabular}
\caption{Parameters of quasi-Least-Favorable Settings
studied in the empirical results presented here.\label{table:formalPred}}
\end{center}
\end{table}

\subsubsection*{Results at $N=1500$}

To test these predictions, we generate in each situation
$R=200$ random realizations of size $N=1500$ from
$\LSF(\delta,\rho,\sigma,\nu)$  with the parameters shown in
Table \ref{table:formalPred}
and run the LARS/LASSO solver to find the solution $\hxl$.
Table \ref{table:Nsmall} shows the empirical average MSE
in $200$ trials at each tested situation.

Except at $\delta = 0.10$ the mismatch between empirical and theoretical
a few to several percent \-- reasonable given the sample size $R=200$.
At $\delta=0.10$, $\rho = 0.180$ -- close to
phase transition --  there is a mismatch needing attention. (In fact, at each level of $\delta$ the most serious mismatch
is at the value of $\rho$ closest to phase transition. This can be attributed partially to the blowup
of the quantity being measured as we approach phase transition.)  We will pursue this mismatch below.

We also ran trials at $\delta \in\{0.15, 0.20,0.30, 0.35, 0.40, 0.45\}$.
These cases exhibited the same patterns seen above, with adequate fit
except at small $\delta$, especially near phase transition. We omit the data here.

In all our trials, we measured numerous observables -- not only the MSE.
The trend in mismatch between theory and observation
 in such observables was comparable to that seen for MSE.
 In \cite{AMPSupplement,TheoPredLassoOpChar}, the reader can find discussion
 and presentation of evidence for other observables.
\begin{table}
\begin{center}
\begin{tabular}{|l|l|l|l|l|l|l|}
\hline
$\delta$ & $\rho$ & $\mu$ & $\lambda^*$ & fMSE & eMSE & SE \\
\hline
0.100 & 0.095 & 5.791 & 1.258 & 0.136 & 0.126 & 0.0029\\
0.100 & 0.142 & 8.242 & 0.804 & 0.380 & 0.329 & 0.0106\\
0.100 & 0.170 & 12.901 & 0.465 & 1.045 & 0.755 & 0.0328\\
0.100 & 0.180 & 18.278 & 0.338 & 2.063 & 1.263 & 0.0860\\
\hline
0.250 & 0.134 & 5.459 & 0.961 & 0.374 & 0.373 & 0.0046\\
0.250 & 0.201 & 7.683 & 0.592 & 1.028 & 1.002 & 0.0170\\
0.250 & 0.241 & 12.219 & 0.351 & 2.830 & 2.927 & 0.0733\\
0.250 & 0.254 & 17.314 & 0.244 & 5.576 & 5.169 & 0.1978\\
\hline
0.500 & 0.193 & 5.194 & 0.689 & 0.853 & 0.836 & 0.0078\\
0.500 & 0.289 & 7.354 & 0.400 & 2.329 & 2.251 & 0.0254\\
0.500 & 0.347 & 11.746 & 0.231 & 6.365 & 6.403 & 0.1157\\
0.500 & 0.366 & 16.667 & 0.159 & 12.427 & 11.580 & 0.2999\\
\hline
\end{tabular}
\caption{Results at $N=1500$. MSE of \LASSO($\lambda^*)$ at
nearly-least-favorable situations, together with standard errors (SE)}
\label{table:Nsmall}
\end{center}
\end{table}

\subsubsection*{Results at $N=4000$}

Statistics of random sampling dictate that there always be some
measure of disagreement between empirical averages and expectations.
When the expectations are taken in the large-system limit, as ours are,
there are additional small-$N$ effects that appear separate from random sampling effects.
However, both sorts of effects should visibly decline with increasing $N$.

Table \ref{table:Nbig} presents results for $N=4000$;
we expect the discrepancies to shrink when the experiments are run at larger value of $N$.
We study the same $\rho$ and $\delta$ that
were studied for $N=1500$, and see that the mismatches
in our MSE's have grown smaller with $N$.
\begin{table}
\begin{center}
\begin{tabular}{|l|l|l|l|l|l|l|}
\hline
$\delta$ & $\rho$ & $\mu$ & $\lambda^*$ & fMSE & eMSE& SE\\
\hline
0.100 & 0.095 & 5.791 & 1.258 & 0.136 & 0.128 & 0.0016\\
0.100 & 0.142 & 8.242 & 0.804 & 0.380 & 0.348 & 0.0064\\
0.100 & 0.170 & 12.901 & 0.465 & 1.045 & 0.950 & 0.0228\\
0.100 & 0.180 & 18.278 & 0.338 & 2.063 & 1.588 & 0.0619\\
\hline
0.250 & 0.134 & 5.459 & 0.961 & 0.374 & .371 & 0.0028\\
0.250 & 0.201 & 7.683 & 0.592 & 1.028 &  1.023 & 0.0106\\
0.250 & 0.241 & 12.219 & 0.351 & 2.830 &  2.703 & 0.0448\\
0.250 & 0.254 & 17.314 & 0.244 & 5.576 &  5.619 & 0.0428\\
\hline
0.500 & 0.193 & 5.194 & 0.689 & 0.853 & 0.849 & 0.0047\\
0.500 & 0.289 & 7.354 & 0.400 & 2.329 & 2.296 & 0.016\\
0.500 & 0.347 & 11.746 & 0.231 & 6.365 & 6.237 & 0.0677 \\
0.500 & 0.366 & 16.667 & 0.159 & 12.427& 12.394 & 0.171 \\
\hline
\end{tabular}
\caption{Results at $N=4000$. Theoretical and empirical MSE's of \LASSO($\lambda^*)$ at
nearly-least-favorable situations, together with standard errors (SE).}
\label{table:Nbig}
\end{center}
\end{table}

\subsubsection*{Results at $N=8000$}
Small values of $\delta$ have the largest discrepancy specially when $\rho$ is chosen very close to the phase transition curve. To show that
this discrepancy shrinks as we increase the value of $N$, we do a similar experiment for $\delta=0.10$ but this time with $N=8000$.
Table \ref{table:Nbigbig} summarizes the results of this simulation and shows better agreement between the formal predictions and empirical results.

\begin{table}
\begin{center}
\begin{tabular}{|l|l|l|l|l|l|l|}
\hline
$\delta$ & $\rho$ & $\mu$ & $\lambda^*$ & fMSE & eMSE & SE\\
\hline
0.100 & 0.095 & 5.791 & 1.258 & 0.136 & 0.131 &  0.0012\\
0.100 & 0.142 & 8.242 & 0.804 & 0.380 & 0.378 &  0.0046\\
0.100 & 0.170 & 12.901 & 0.465 & 1.045 & 1.024 & 0.0186\\
0.100 & 0.180 & 18.278 & 0.338 & 2.063 & 1.883 &  0.0458\\
\hline
\end{tabular}
\caption{Results at $N=8000$. Theoretical and empirical MSE's of \LASSO($\lambda^*)$ at
nearly-least-favorable situations with $\delta=0.10$,
together with standard errors (SE) of the empirical MSE's}
\label{table:Nbigbig}
\end{center}
\end{table}

\begin{figure}
\center{\includegraphics[width=10.cm]{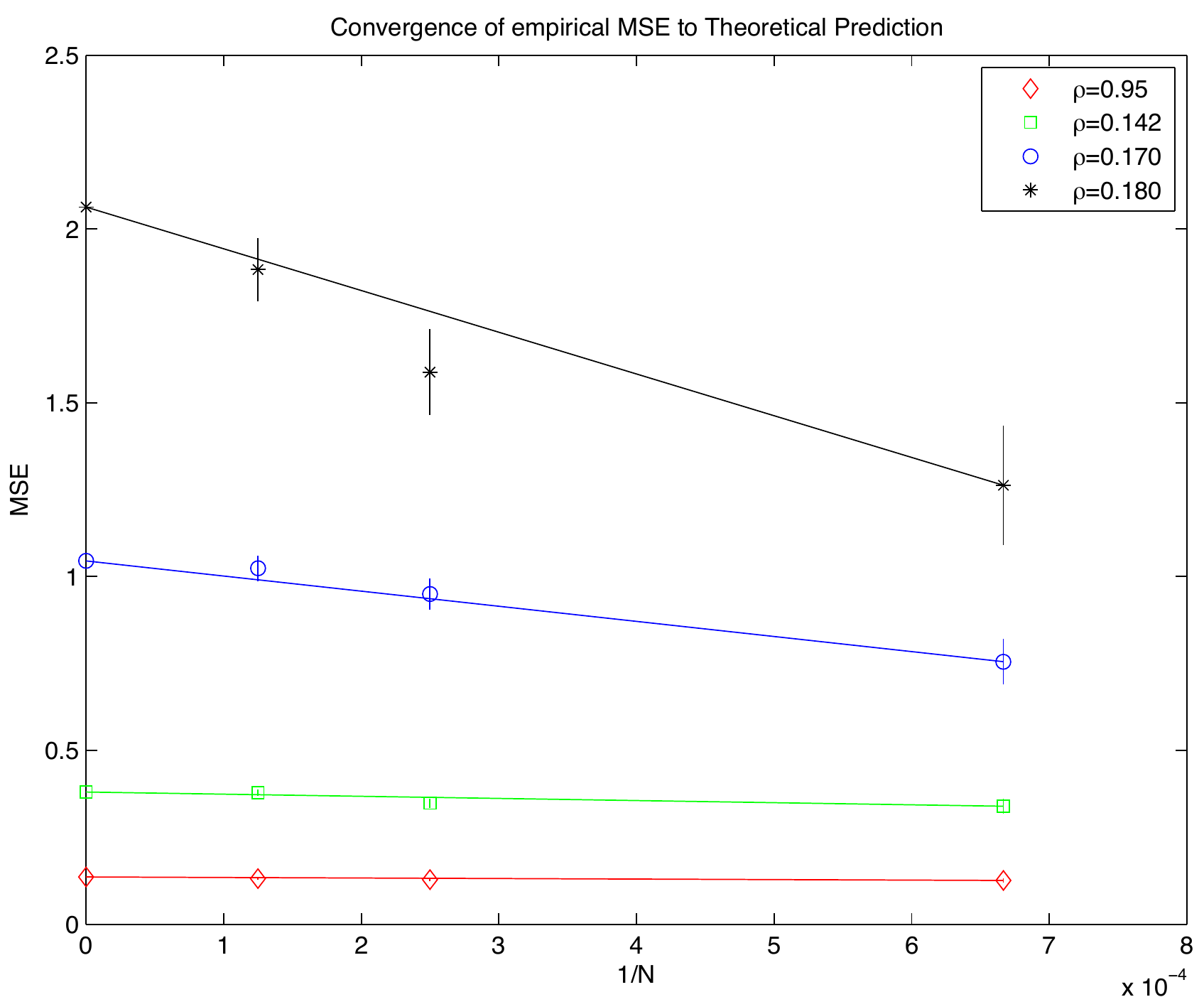}}
\caption{Finite-$N$ scaling of empirical MSE.  Empirical MSE results
from the cases $N=1500$, $N=4000$ and $N = 8000$ and $\delta=0.1$.
Vertical axis: empirical MSE. Horizontal axis: $1/N$. Different colors/symbols
indicate different values of the sparsity control parameter $\delta$. Vertical bars denote $\pm 2  SE$
limits. Theoretical predictions for the $N = \infty$ case appear at $1/N = 0$. Lines connect
the cases $N = 1500$ and $N=\infty$. }
\label{fig:Finite-N-scaling}
\end{figure}

The alert reader will no doubt have noticed that the discrepancy between theoretical predictions
and empirical results is in many cases quite a bit larger in magnitude than the
size of the the formal standard errors reported in the above tables.  We emphasize
that the theoretical predictions are formal limits for the $N \goto \infty$ case, while
empirical results take place at finite $N$.  In both statistics and statistical physics it is
quite common for mismatches between finite-$N$ results and $N$-large to occur as
either $O(N^{-1/2})$ (eg Normal approximation to the Poisson) or $O(N^{-1})$ effects
(eg Normal approximation to fair coin tossing).  Analogously, we might anticipate that
 mismatches in this setting of order $N^{-\alpha}$ with $\alpha$ either $1/2$ or $1$.
Figure \ref{fig:Finite-N-scaling} presents empirical and theoretical results
taken from the cases $N=1500$, $4000$, and $8000$ and displays them on  a common
graph, with $y$-axis  a mean-squared error (empirical or theoretical) and on the $x$ axis
the inverse system size $1/N$.  The case $1/N = 0$ presents the formal large-system
limit predicted by our calculations and the other cases $1/N > 0$ present empirical results
described in the tables above. As can be seen, the discrepancy between formal MSE
and empirical MSE tends to zero  linearly with $1/N$.  (A similar plot with $1/\sqrt{N}$ on
the $x$-axis would not be so convincing.)

\begin{finding}   The formal and empirical $\MSE$'s  at the quasi saddlepoint $(\nu^*,\lambda^*)$
show statistical agreement
at the cases studied, in the sense that either the $\MSE$'s are consistent with standard
statistical sampling formulas, or, where they were not consistent
at $N=1500$, fresh data at $N=4000$ and $N=8000$ showed marked reductions in the anomalies
confirming that the anomalies decline with increasing $N$.
\end{finding}

\subsubsection{Existence of Game-Theoretic Saddlepoint in eMSE}

\begin{figure}
\center{\includegraphics[width=10.cm]{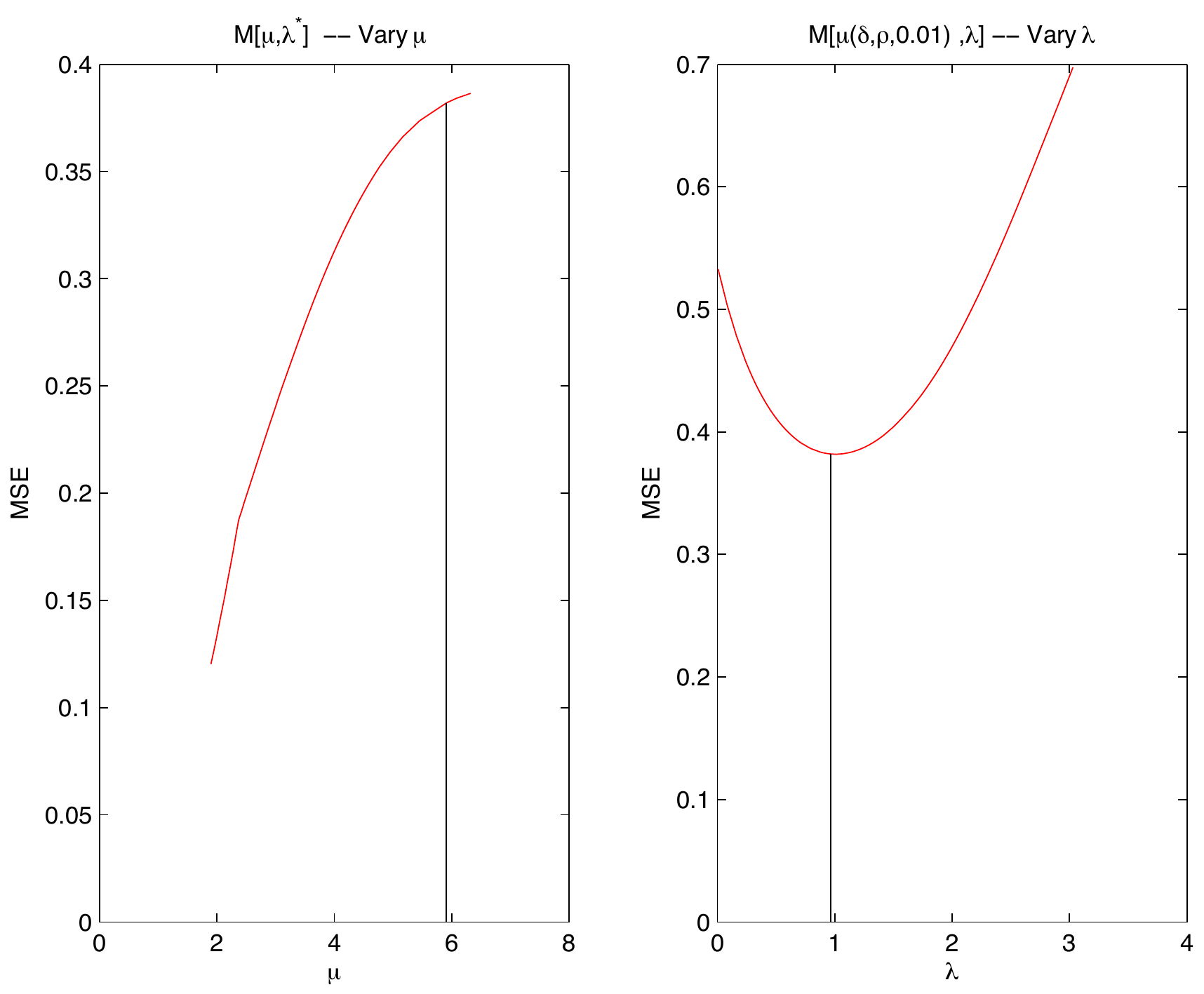}}
\caption{Saddlepoint in formal MSE. Right panel: Behavior of formal MSE
as $\lambda$ is varied away from $\lambda^*$. Left panel: Behavior of
formal MSE as $\mu$ is varied away from $\mu^*$ in the direction of smaller values. Black lines indicate
locations of $\mu^*$ and $\lambda^*$. $\delta=0.25$, $\rho = \rhoMSE(\delta)/2$.}
\label{fig:Saddlepoint}
\end{figure}

Underlying our derivations of minimax formal MSE
is a game-theoretic saddlepoint structure, illustrated in
Figure \ref{fig:Saddlepoint}.  The loss function MSE has the following structure
around the quasi saddlepoint $(\nu^*,\lambda^*)$: any variation
of $\mu$ to lower values, will cause a reduction in loss,
while a variation of $\lambda$ to other values will cause an increase in loss.

\subsubsection{Other penalization gives larger MSE}

If our formalism is correct in deriving optimal penalization for $\hxl$,
we will see that changes of the penalization away from $\lambda^*$ will cause
MSE to increase.
We consider the same situations as earlier, but now vary $\lambda$ away from the
minimax value, while holding the other aspects of the problem fixed.
In the Appendix, Tables  \ref{table:FullVaryLambdaOneFiveK}
and \ref{table:FullVaryLambdaFourK} presents numerical values of the
empirical MSE obtained.
 Note the agreement of formal MSE,
in which a saddlepoint is rigorously proven, and empirical MSE,
which represents actual LARS/LASSO reconstructions. Also in this case
we used $R=200$ Monte Carlo replications.

\begin{figure}
\center{\includegraphics[width=10.cm]{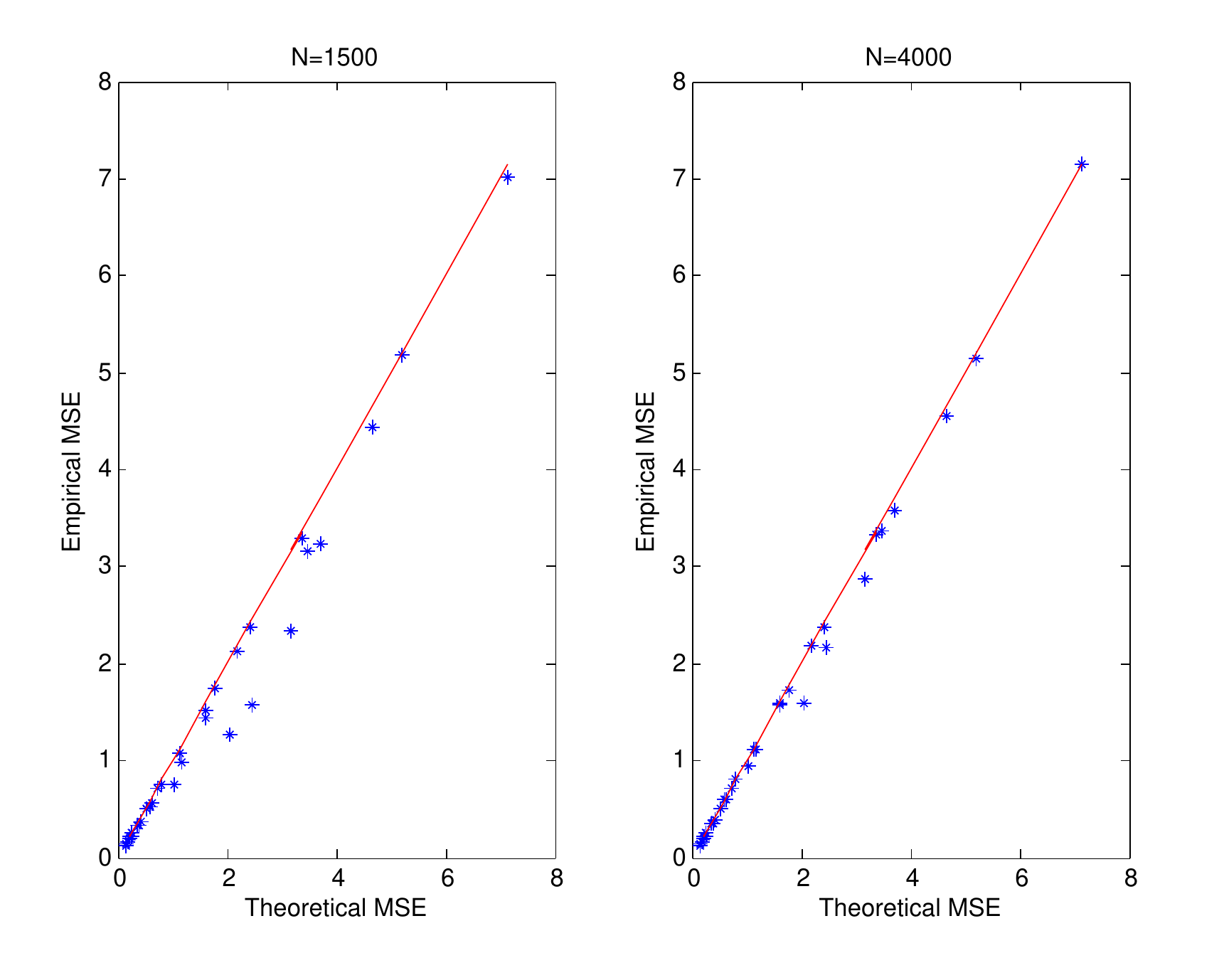}}
\caption{Scatterplots comparing Theoretical and Empirical MSE's found in Tables
 \ref{table:FullVaryLambdaOneFiveK}  and \ref{table:FullVaryLambdaFourK}.
 Left Panel: results at $N=1500$.
 Right Panel: results at $N=4000$.
 Note visible tightening of the scatter around the identity line as $N$ increases.}
\label{fig:ScatterVaryLambda}
\end{figure}

  To visualize the information in those tables, we refer to
Figure \ref{fig:ScatterVaryLambda}.

\subsubsection{MSE with more favorable measures is smaller}

In our formalism, fixing  $\lambda = \lambda^*$,  and varying
$\mu$ to smaller values
will cause a reduction in formal MSE. Namely, if instead of
$\mu^*(\delta,\rho,0.01)$ we
used $\mu^*(\delta,\rho,\alpha)$ for $\alpha$ significantly
larger than $0.01$,
we would see a significant reduction in MSE, by an amount matching the predicted amount.

Recall that $\stMSE(\nu,\tau)$ denotes the `risk' (MSE) of scalar soft
thresholding as in Section \ref{sec:Scalar}, with input distribution
$\nu$, noise variance $1$, and threshold $\tau$. Now suppose that
$\stMSE(\nu_0,\tau)  > \stMSE(\nu_1,\tau)$.
Then also the resulting formal noise-plus-interference obeys
$\FNPI(\nu_0,\tau) > \FNPI(\nu_1, \tau)$.
As noticed several times in Section \ref{sec:MainPropo},
the formal MSE of AMPT obeys $\FMSE(\nu,\tau)  =
\stMSE(\tilde{\nu},\tau)  \cdot \FNPI(\nu,\tau)$, where $\tilde{\nu}$ denotes
a rescaled probability measure
(as in the proof of Proposition \ref{prop:MMxFormula}).
Hence
\[
\FMSE(\nu_1,\tau) \leq \stMSE(\tilde{\nu_1},\tau) \cdot
\FNPI(\nu_0,\tau)\, ,
\]
where the scaling uses $\FNPI(\nu_0)$.   In particular, for $\mu =
\mu^*(\delta,\rho,\alpha) = \mu^{\hash}(\delta\cdot\rho,\alpha)
\sqrt{\NPI^*(\delta,\rho)}$,
the three point mixture: $\nu_{\delta,\rho,\alpha} $ has
\[
    \FMSE(\nu_{\delta,\rho,\alpha},\tau^*) \leq (1-\alpha) M^*(\delta, \rho),
\]
and we ought to be able to see this.  Table \ref{table:VaryMu} shows results of simulations at $N=1500$.
The theoretical MSE drops as we move away from the nearly least favorable $\mu$ in the direction
of smaller $\mu$, and the empirical MSE responds similarly.

\begin{finding}  The empirical data exhibit the saddlepoint structures
predicted by the SE formalism.
\end{finding}
%
%
\subsubsection{MSE of Mixtures}

The SE formalism contains a basic mathematical structure
which allows one to infer that behavior at one saddlepoint
determines the global minimax value: behavior under
taking convex combinations (mixtures) of measures $\nu$.

Let $\stMSE(\nu,\lambda)$ denote the `risk' (MSE) of scalar soft thresholding
as in Section 2. For such scalar thresholding, we have the affine relation
\[
\stMSE((1-\gamma) \nu_0 + \gamma \nu_1,\tau) = (1-\gamma)\stMSE(\nu_0,\tau) + \gamma \cdot \stMSE(\nu_1,\tau)\, .
\]
Now suppose that  $\stMSE(\nu_0,\tau)  > \stMSE(\nu_1,\tau)$.
Then also $\NPI(\nu_0,\tau) > \NPI(\nu_1, \tau)$.
The formal MSE of AMPT obeys the scaling
relation $\FMSE(\nu,\tau)  = \stMSE(\tilde{\nu},\tau)  \cdot
\NPI(\nu,\tau)$, where $\tilde{\nu}$ denotes
the rescaled probability measure, argument rescaled by $1/\sqrt{NPI}$.
We conclude that
\begin{equation} \label{eq:convexitybound}
\FMSE((1-\gamma) \nu_0 + \gamma \nu_1,\tau)  \leq (1-\gamma) \cdot
\stMSE(\tilde{\nu}_0,\tau) \cdot \NPI(\nu_0,\tau)+ \gamma \cdot \stMSE(\tilde{\nu}_1,\tau) \cdot \NPI(\nu_0,\tau),
\end{equation}
This  `quasi-affinity' relation allows to extend
the saddlepoint structure
from 3 point mixtures to more general
measures.

\begin{figure}
\center{\includegraphics[width=10.cm]{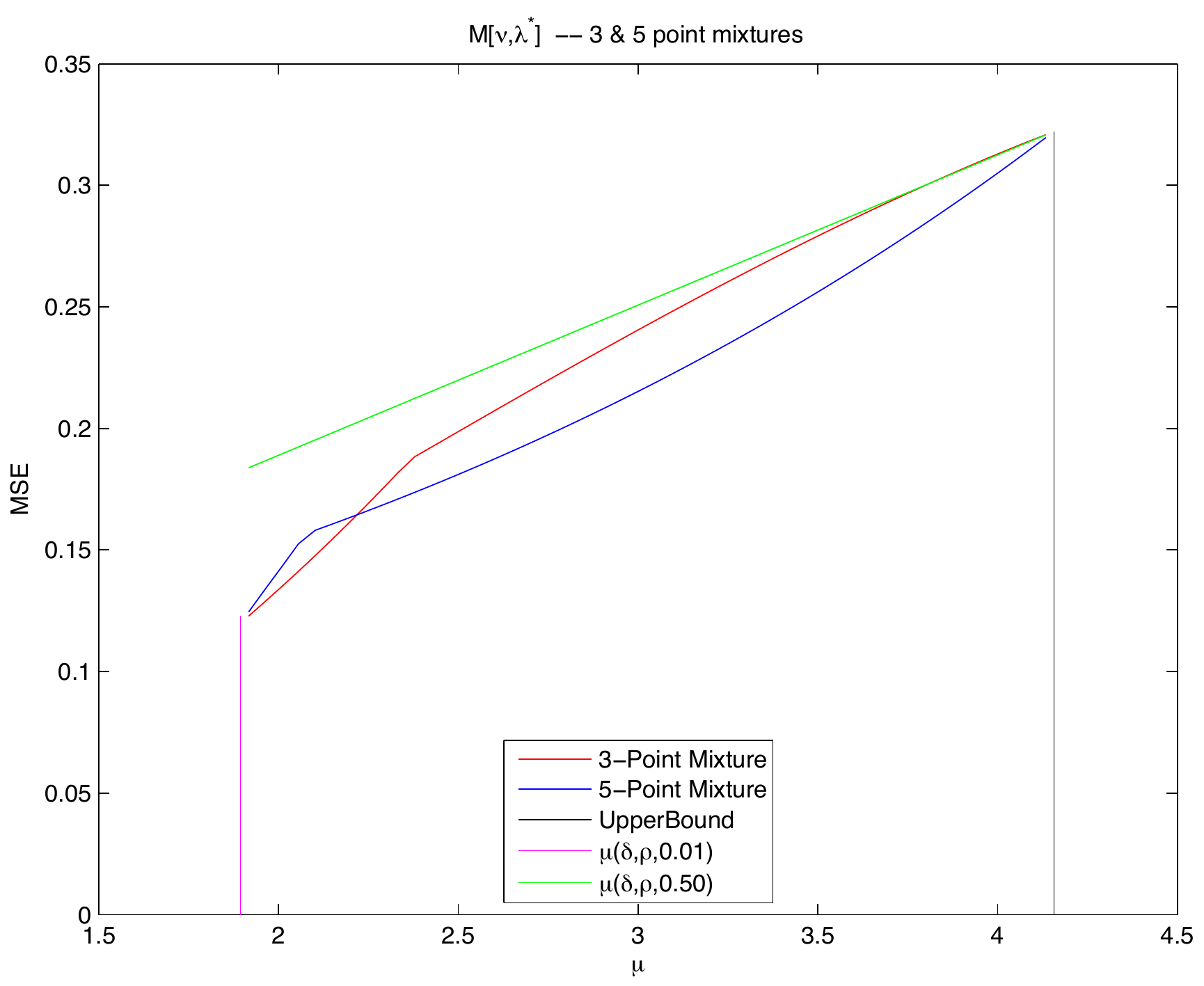}}
\caption{Convexity structures in formal MSE. Behavior of formal MSE
of 5 point mixture combining nearly least-favorable $\mu$ with discount of 1\%
and one with discount of 50\%. Also, the convexity bound (\ref{eq:convexitybound})
and the formal MSE of associated 3-point mixtures is displayed.
$\delta=0.25$, $\rho = \rhoMSE(\delta)/2$.}
\label{fig:Convexity}
\end{figure}

To check this, we consider two near-least-favorable measures, $\nu_0 = \nu_{\delta,\rho,0.02}$ and $\nu_1 = \nu_{\delta,\rho,0.50}$.
and generate a range of cases $\nu^{(\alpha)} = (1-\alpha) \nu_0 + \alpha \nu_1$ by varying alpha.
When $\alpha \not\in \{0,1\}$ this is a 5 point mixture rather than one of the 3-point mixtures we have been studying.
Figure \ref{fig:Convexity} displays the convexity bound (\ref{eq:convexitybound}),
and the behavior of the formal MSE of this 5 point mixture. For comparison it also presents the formal MSE
of the 3 point mixture having its mass at the weighted mean  $(1-\alpha) \mu(\delta,\rho,0.02) + \alpha \mu(\delta,\rho,0.50)$.
Evidently, the 5 point mixture typically has smaller MSE than the comparable 3-point mixture, and it always
is below the convexity bound.

\begin{finding}  The empirical MSE obeys the mixture inequalities
predicted by the SE formalism.
\end{finding}

\subsection{Above Phase Transition}

We conducted an empirical study of the formulas derived in Section  4.5.
At  $\delta = 0.25$ we chose $\rho = 0.401$ \-- well above phase transition \--
and selected a range of $\tau$ and $\gamma$ values allowed by
our formalism.  For each pair $\gamma,\tau$, we generated $R=200$
Monte Carlo realizations and obtained LASSO solutions with the
given penalization parameter $\lambda$.  The results are
described in Table \ref{table:AbovePT}.  The match between
formal MSE and empirical MSE is acceptable.

\begin{table}[h]
\begin{center}
\begin{tabular}{|r|r|r|r|r|r|r|r|}
\hline
 $\delta$ & $\rho$ &  $\gamma$ & $\mu$ & $\tau$ & $\lambda$ & fMSE & eMSE \\
\hline
0.250 & 0.401 & 0.75& 2.8740  & 1.500 &  0.9840  & 0.750 & 0.746  \\
0.250 & 0.401  & 0.85 & 4.142 & 1.500 &  1.168  & 1.417 & 1.425\\
0.250 & 0.401 & 0.90 & 5.345 & 1.500 & 1.366 & 2.250 & 2.239 \\
0.250 & 0.401 &  0.95 & 7.954 & 1.500 & 1.841 & 4.750 & 4.724 \\
0.250 & 0.401 &  0.97 & 10.4781 & 1.500 &  2.328 & 8.083 & 8.126 \\
0.250 & 0.401 &  0.98 & 12.9628  & 1.500 &  2.822 & 12.250 &12.327 \\
0.250 & 0.401 &  0.99& 18.5172 & 1.500 &  3.949 & 24.750 & 24.601\\
0.250 & 0.401 &  0.995 & 26.3191 & 1.500 &  5.5558 &49.750 & 49.837\\
\hline
0.250 & 0.401 &0.75 & 2.9031 & 2.000 &  2.8766 & 1.417 & 1.409 \\
0.250 & 0.401 &  0.85& 4.058 & 2.000 & 3.626 & 2.250 & 2.238 \\
0.250 & 0.401 &  0.90 & 5.158 & 2.000 & 4.385 & 2.250 & 2.238 \\
0.250 & 0.401 & 0.95 & 7.560 & 2.000 & 6.122 & 4.750 & 4.742 \\
0.250 &0.401   & 0.97 & 9.897 & 2.000 & 7.861 & 8.083 & 8.054   \\
0.250 & 0.401  & 0.98 & 12.205 & 2.000 & 9.6019 & 12.250 & 12.215 \\
0.250 & 0.401  & 0.99 & 17.380 & 2.000 & 13.5425 & 24.750 & 24.634 \\
0.250 & 0.401  & 0.995& 24.662 & 2.000 & 19.1260 & 49.750  & 49.424 \\
\hline
0.250 & 0.401 & 0.75  & 2.817 & 2.500 & 4.501 & 1.417 & 1.409 \\
0.250 & 0.401 &  0.85 & 3.896 & 2.500 & 5.750 & 2.250 & 2.241 \\
0.250 & 0.401 &  0.90 & 4.926 & 2.500 & 7.004 & 2.250 & 2.241 \\
0.250 & 0.401 & 0.95 & 7.181 & 2.500 & 9.848 & 4.750 & 4.712 \\
0.250 & 0.401 &  0.97  & 9.380 & 2.500 & 12.6846 & 8.083 & 8.050   \\
0.250 & 0.401 &  0.98  & 11.555 & 2.500 & 15.5170 & 12.250 & 12.215  \\
0.250 & 0.401 & 0.99  & 16.436 & 2.500 & 21.9183 & 24.750 & 24.619   \\
0.250 & 0.401 &0.995 &  23.311 & 2.500 &30.9786 & 49.750 & 49.442 \\
\hline
0.250 & 0.401 & 0.75 &  2.7649 & 3.000 & 5.8144& 1.417 & 1.408 \\
0.250 & 0.401 &  0.85 & 3.809 & 3.000 & 7.4730 & 2.250 & 2.241 \\
0.250 & 0.401 &  0.90 & 4.806 & 3.000 & 9.131 & 2.250 & 2.241 \\
0.250 & 0.401 & 0.95  & 6.991 & 3.000 & 12.880 & 4.750 & 4.735 \\
0.250 & 0.401 & 0.97 & 9.125  & 3.000   & 16.6113 & 8.083      & 8.053 \\
0.250 & 0.401 &0.98 & 11.236  & 3.000  & 20.3339 & 12.250 & 12.218 \\
0.250 & 0.401 &0.99 & 15.975   & 3.000 & 28.7413 & 24.750 & 24.621\\
0.250 & 0.401 &0.995 & 22.652  & 3.000 & 40.6356 & 49.750 & 49.419 \\
 \hline
\end{tabular}
\caption{Results above Phase transition.  Parameters of the construction as well
as theoretical predictions and resulting empirical MSE figures}
\label{table:AbovePT}
\end{center}
\end{table}

\begin{finding}
Running $\hxl$ at the $3$-point mixtures defined
for the regime above phase transition
in Lemma \ref{lem:LASSOAbovePT}
yields empirical MSE consistent with the formulas of that Lemma.
\end{finding}
This validates the unboundedness of MSE of LASSO above phase transition.
%
%

\section{Extensions}

\subsection{Positivity Constraints}
\label{sec-Positivity}

A completely parallel treatment can be given for the case where $x^0 \geq 0$.
In that setting, we use the positivity-constrained soft-threshold
\begin{eqnarray}
\eta^+(x;\theta) = \left\{\begin{array}{ll}
x-\theta & \mbox{ if $\theta<x$,}\\
0     & \mbox{ if $ x\le \theta$,}\\
\end{array}\right.
\end{eqnarray}
and consider the corresponding
positive-constrained thresholding minimax MSE \cite{DJHS92}
\begin{eqnarray}
    M^+(\eps) = \inf_{\tau > 0 } \sup_{\nu \in \cF^+_\eps}
 \E\Big\{\big[\eta^+\big(X+ {\sigma} \cdot Z;
\tau\sigma\big)-X\big]^2\Big\} ,
\end{eqnarray}
where
\[
   \cF^+_{\eps}  = \{ \nu\,:\, \nu \mbox{ is probability measure with }
 \nu[0,\infty) = 1, \nu(\{0\}) \ge 1-\eps \}.
\]
We consider the positive-constrained $\ell_1$-penalized least-squares
estimator $x^{1,\lambda,+}$, the solution to
\beq \label{Pos}
(P_{2,\lambda,1}^+) \qquad  {\rm minimize}_{x \geq 0} \quad  \frac{1}{2}\, \| y - Ax\|_2^2 +
\lambda \|x\|_1 .
 \eeq
We define the minimax, formal {\it noise sensitivity}:
\begin{equation} \label{eq:DefMStarPlus}
     M^{+,*}(\delta,\rho) = \sup_{\sigma > 0}  \max_\nu    \min_\lambda \FMSE(x^{1,\lambda,+},\nu,\sigma^2)/\sigma^2;
\end{equation}
here $\nu \in \cF^+_{\rho\delta}$ is the marginal distribution of $x_0$.
Let $\rhoMSE^+(\delta)$ denote the solution of
 \begin{eqnarray} \label{DefRhoPlusMSE}
M^+(\rho\delta) = \delta\, .
\end{eqnarray}
In complete analogy to (\ref{eqMStarResult}) we have the formula:
\beq
   M^{+,*}(\delta,\rho) =  \left\{ \begin{array}{ll}
                                         \frac{M^+(\delta\rho)}{1-M^+(\delta\rho)/\delta} , & \rho < \rhoMSE^+(\delta) ,\\
                                         \infty  ,& \rho \geq  \rhoMSE^+(\delta). \\
                                           \end{array} \right.
\eeq

\begin{figure}
\center{\includegraphics[width=10.cm]{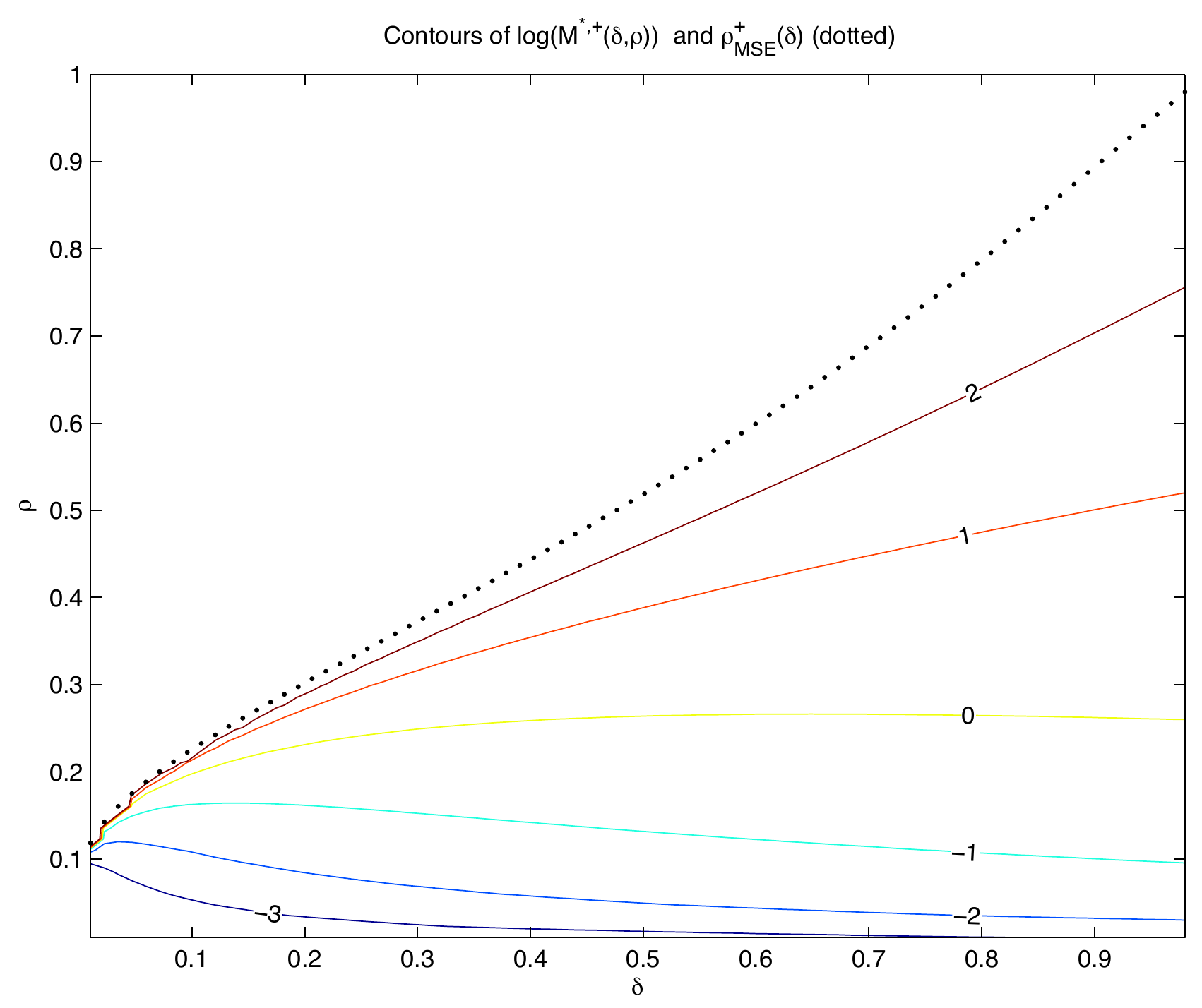}}
\caption{{\small Contour lines of the positivity-constrained minimax noise sensitivity
$\Ms^{*,+}(\delta, \rho)$ in the $(\rho,\delta)$ plane. The
dotted black curve graphs the phase boundary $(\delta,\rhoMSE^+(\delta))$.
Above this curve, $\Ms^{*,+}(\delta, \rho) =\infty$. The
colored lines present level sets of $\Ms^{*,+}(\delta, \rho)
=1/8$, $1/4$, $1/2$, $1$, $2$, $4$ (from bottom to top).}}
\label{fig:Figure001Pos}
\end{figure}

The argument is the same as above, using the AMP formalism,
with obvious modifications. The papers \cite{DMM,AMPSupplement} show in more detail
how to make arguments for AMP that apply simultaneously to the sign-constrained
and unconstrained case. All other features of Proposition \ref{prop:main}
carry over, with obvious substitutions.  Figure \ref{fig:Figure001Pos} shows the
phase transition for the positivity constrained case, as well as the contour lines of $M^{+,*}$.
Again in analogy to the sign-unconstrained case, the phase boundary $\rho_{MSE}^+$
occurs at precisely the same location at the phase boundary for $\ell_1$-$\ell_0$ equivalence;
as earlier this can be inferred from formulas in this paper and in \cite{DMM}.

\subsection{Other Classes of Matrices}

We focused here on matrices $A$ with Gaussian iid entries.

Previously, extensive empirical evidence was presented by Donoho and Tanner \cite{DoTa08_universality},
that pure $\ell_1$-minimization has its $\ell_1$-$\ell_0$ equivalence phase transition
at the boundary $\rhoMSE^{\pm}$ not only for Gaussian matrices but for a wide
collection of ensembles, including partial Fourier, partial Hadamard, expander graphs,
iid $\pm 1$. This is the noiseless, $\lambda=0$ case of the general noisy, $\lambda \geq 0$ case studied here.

We believe that similar results to those obtained here
hold for matrices $A$ with uniformly bounded iid entries with zero
mean and variance $1/n$. In fact, we believe our results should extend to a broader
universality class including matrices with iid entries with same mean and
variance, under an appropriate light tail condition.

\section{Relations with Statistical Physics and Information Theory}
\label{sec:StatPhys}

This section outlines the relations of the approach advocated
here with ideas in information theory (in particular, with the theory
of sparse graph codes), graphical models and statistical physics
(more precisely spin glass theory). We will not discuss such
relations in full mathematical detail, but only stress some important
points that might be useful for researchers in each of those fields.
%
%
\subsection{Information theory and message passing algorithms}

Message passing algorithms, and most notably belief propagation,
have been intensively investigated in coding theory and communications,
in particular because of their success in decoding sparse graph
codes \cite{RiUBook}. Belief propagation is defined whenever the
\emph{a posteriori} joint distribution of the variables to be inferred
$x$ conditional on the observations $y$ can be written
as a graphical model.
In the present case this is easily done, provided the a priori
probability distribution of the signal
$x=(x_1,\dots,x_N)$ takes a
 $\nu = \nu_1\times \nu_2 \dots \times \nu_N$.
The posterior is then
\begin{eqnarray*}
\mu(\de x)= \frac{1}{Z}  \prod_{a=1}^n \exp\Big\{-\frac{\beta}{2}(y_a-
(Ax)_a)^2\Big\}
\prod_{i=1}^N \nu_i(\de x_i)\, .\label{eq:Posterior}
\end{eqnarray*}
Graphical models of this type were (implicitly or explicitly)
considered in the context of  multiuser detection
\cite{Kabashima,NeirottiSaad,MoTse1,MoTse2}. The underlying factor graph
\cite{KFL01} is the complete bipartite graph over $N$ variable nodes
and $n$ factor nodes.

Applying belief propagation to such a model incurs two
obvious difficulties: the graph is dense (and hence the complexity
per iteration scales at least like $n^3$, and in fact worse), and the
alphabet is continuous (and hence messages are not finitely representable).
As discussed in \cite{DMM_ITW1}, AMP solves these problem.
From the information theoretical perspective,
the term $+z^{t-1} \df_{t}/n$ in Eq.~(\ref{eq:FOAMP2})
corresponds to `subtracting intrinsic information'.

An important difference between the message passing algorithms in coding theory and
what is presented here is that no precise information is available
on the priors $\nu_i$ in Eq.~(\ref{eq:Posterior}). Therefore
the AMP rules should not be sensitive to the prior. The use of the soft
threshold function $\eta(\,\cdot\,;\theta)$ makes the AMP robust
within the class of sparse priors. Also, it is directly related to
the $\ell_1$ regularization in the LASSO.

In coding theory, message passing algorithms
are analyzed through density evolution \cite{RiUBook}.
The common justification for density evolution is that the underlying
graph is random and sparse, and hence converges locally to a tree
in the large system limit. In the case of trees density
evolution is exact, hence it is asymptotically exact for sparse
random graphs. Such an easy justification is not available
in the cases of dense graphs treated here and a deeper
mathematical analysis is required.
In \cite{BayatiMontanari}, this analysis was carried out
in the case of Gaussian matrices $A$. It remains a challenge
to generalize such analysis beyond the case of Gaussian matrices $A$.

Having outlined the relation with belief propagation and coding,
it is important to clarify a key point. In the context of
sparse graph coding, belief propagation performances
and MAP (maximum a posteriori probability) performances
do not generally coincide even asymptotically (although
they are intimately related \cite{MMU05,MMRU05}).
In the present paper we instead conjecture that AMP and LASSO
have asymptotically equal MSE under appropriate calibration. This is
due to the fact that that the state evolution recursion
$m_t\mapsto m_{t+1}=\Psi(m_t)$ has only one stable fixed point.
%
%
\subsection{Statistical physics}

There is a well studied connection between statistical physics
techniques and message passing algorithms \cite{MezardMontanari}.
In particular, the sum-product algorithm corresponds to
the Bethe-Peierls approximation  in statistical physics, and its
fixed points are stationary points of the Bethe free energy.
In the context of  spin glass theory, the
Bethe-Peierls approximation is also referred to as the
`replica symmetric cavity\footnote{When this terminology is used in
statistical physics,
the emphasis is rather on properties of random instances.}
method'.

The Bethe-Peierls approximation postulates a set
of non-linear equations on quantities that correspond to the
belief propagation messages, and allow to compute posterior
marginals under the distribution (\ref{eq:Posterior}).
In the special cases of
spin glasses on the complete graph (the celebrated Sherrington-Kirkpatrick
model), these equations reduce to the so-called TAP equations,
named after Thouless, Anderson and Palmer who first used them
\cite{TAP}.

The original TAP equations where a set of non-linear
equations for local magnetizations (i.e. expectations of a single variable).
Thouless, Anderson and Palmer first recognized that naive mean field
is not accurate enough in the spin glass model, and corrected it
by adding the so called Onsager reaction term that is analogous
to the term  $+z^{t-1} \df_{t}/n$
in Eq.~(\ref{eq:FOAMP2}).
More than 30 years after the original paper, a complete mathematical
justification of the TAP equations remains an
open problem in spin glass theory,
although important partial results exist \cite{Talagrand}.
While the connection between belief propagation and
Bethe-Peierls approximation stimulated a considerable
amount of research \cite{YedidiaFreemanWeiss},
the algorithmic uses of TAP equations
have received only sparse attention. Remarkable exceptions
include \cite{OpperWinther,Kabashima,NeirottiSaad}.
%
%
\subsection{State evolution and replica calculations}

Within statistical mechanics, the typical properties of
probability measures of the form (\ref{eq:Posterior})
are studied using the replica method or the cavity method
\cite{MezardMontanari}. These can be described as non-rigorous but
mathematically sophisticated techniques. Despite intense efforts
and some spectacular progresses \cite{Talagrand}, even a
precise statement of the assumptions implicit in such techniques
is missing, in a general setting.

The fixed points of state evolution describe the
output of the corresponding AMP, when the latter
is run for a sufficiently large number of iterations
(independent of the dimensions $n,N$).
It is well known, within statistical mechanics \cite{MezardMontanari},
that the fixed point equations do indeed coincide with
the equations obtained form the replica method
(in its replica-symmetric form).

During the last few months, several papers investigated
compressed sensing problems using the replica method
\cite{Goyal,KabashimaTanaka,BaronGuoShamai}.
In view of the discussion above, it is not surprising
that these results can be recovered from the state evolution formalism
put forward in \cite{DMM}. Let us mention that the latter
has several advantages over the replica method:
\begin{itemize}
\item[$(1)$] It is more concrete, and its assumptions can be
checked quantitatively through simulations;
\item[$(2)$] It is intimately related to efficient message passing algorithms;
\item[$(3)$] It actually allows to predict the performances of these
algorithms (including for instance precise convergence
time estimates);
\item[$(4)$] It actually leads to rigorous
statements, at least in the case of Gaussian sensing matrices.
\end{itemize}

\appendix
%
%
\section{Some explicit formulae}
\label{app:Formulae}

This appendix contain some formulae and analytical derivation
omitted from the main text.

The phase boundary curve admits the parametric expression
\begin{eqnarray}
\delta & = &\frac{2\phi(\tau)}{\tau+2(\phi(\tau)-\tau\Phi(-\tau))}\, ,\\
\rho & = &1-\frac{\tau\Phi(-\tau)}{\phi(\tau)}\, ,\\
\end{eqnarray}
This is simply obtained from Eq.~(\ref{eq:SupExpression}). If we call
$G_\eps(\tau)$ the function on the right hand side, then the parametric
expression given here follows from
$\delta = G_{\eps}(\tau)$ and $G'_{\eps}(\tau)=0$ (which
are equivalent to $\delta = M^\hash(\eps)$).

%
%
\section{Tables}

This appendix contains table of empirical results supporting
our claims.

\begin{table}[t]
\phantom{a}\hspace{0cm}

\vspace{-3.5cm}
\caption{$N=1500$, $\lambda$ dependence of the $\MSE$ at fixed $\mu$
\label{table:FullVaryLambdaOneFiveK}}
\vspace{0.5cm}
\centering 
{\small
\begin{tabular}{|l|l|l|l|l|l|l|l|l|l|l|l|} 
\hline 
$\delta$& $\rho$ & $\mu$ & $\lambda$ & fMSE & eMSE & SE\\
 [1ex] 
\hline 
  0.100 & 0.095 & 5.791 & 0.402 & 0.152 & 0.140 & 0.0029\\
 0.100 & 0.095 & 5.791 & 1.258 & 0.136 & 0.126 & 0.0029\\
 0.100 & 0.095 & 5.791 & 2.037 & 0.142 & 0.133 & 0.0030\\
 0.100 & 0.095 & 5.791 & 3.169 & 0.174 & 0.164 & 0.0028\\
 0.100 & 0.095 & 5.791 & 4.948 & 0.239 & 0.228 & 0.0025\\
 \hline
 0.100 & 0.142 & 8.242 & 0.804 & 0.380 & 0.329 & 0.0106\\
 0.100 & 0.142 & 8.242 & 1.960 & 0.408 & 0.374 & 0.0087\\
 0.100 & 0.142 & 8.242 & 3.824 & 0.534 & 0.504 & 0.0084\\
 0.100 & 0.142 & 8.242 & 6.865 & 0.737 & 0.716 & 0.0059\\
 \hline
 0.100 & 0.170 & 12.906 & 0.465 & 1.045 & 0.755 & 0.0328\\
 0.100 & 0.170 & 12.906 & 2.298 & 1.178 & 0.992 & 0.0326\\
 0.100 & 0.170 & 12.906 & 5.461 & 1.619 & 1.520 & 0.0273\\
 0.100 & 0.170 & 12.906 & 10.607 & 2.197 & 2.138 & 0.0139\\
 \hline
 0.100 & 0.180 & 18.278 & 0.338 & 2.063 & 1.263 & 0.0860\\
  0.100 & 0.180 & 18.278 & 2.934 & 2.467 & 1.573 & 0.0741\\
 0.100 & 0.180 & 18.278 & 7.545 & 3.474 & 3.167 & 0.0569\\
 0.100 & 0.180 & 18.278 & 14.997 & 4.677 & 4.438 & 0.0321\\
\hline
\hline
 0.250 & 0.134 & 5.459 & 0.518 & 0.403 & 0.390 & 0.0044\\
 0.250 & 0.134 & 5.459 & 0.961 & 0.374 & 0.373 & 0.0046\\
 0.250 & 0.134 & 5.459 & 1.419 & 0.385 & 0.386 & 0.0046\\
 0.250 & 0.134 & 5.459 & 2.165 & 0.452 & 0.455 & 0.0053\\
 0.250 & 0.134 & 5.459 & 3.555 & 0.623 & 0.612 & 0.0042\\
 \hline
 0.250 & 0.201 & 7.683 & 0.036 & 1.151 & 1.155 & 0.0174\\
 0.250 & 0.201 & 7.683 & 0.592 & 1.028 & 1.002 & 0.0170\\
 0.250 & 0.201 & 7.683 & 1.183 & 1.073 & 1.069 & 0.0169\\
 0.250 & 0.201 & 7.683 & 2.243 & 1.324 & 1.293 & 0.0158\\
 0.250 & 0.201 & 7.683 & 4.392 & 1.861 & 1.837 & 0.0114\\
 \hline
 0.250 & 0.241 & 12.219 & 0.351 & 2.830 & 2.927 &0.0733\\
 0.250 & 0.241 & 12.219 & 1.219 & 3.065 & 2.998 &0.0661\\
 0.250 & 0.241 & 12.219 & 2.917 & 4.055 & 4.020 & 0.0485\\
 0.250 & 0.241 & 12.219 & 6.444 & 5.709 & 5.625 & 0.0330\\
 \hline
 0.250 & 0.254 & 17.314 & 0.244 & 5.576 & 5.169 & 0.1978\\
 0.250 & 0.254 & 17.314 & 1.433 & 6.291 & 5.992 & 0.1712\\
 0.250 & 0.254 & 17.314 & 3.855 & 8.667 & 8.492 & 0.1148\\
 0.250 & 0.254 & 17.314 & 8.886 & 12.154 & 11.978 & 0.0697\\
\hline \hline
 0.500 & 0.193 & 5.194 & 0.176 & 1.121 & 1.108 & 0.0080\\
 0.500 & 0.193 & 5.194 & 0.470 & 0.894 & 0.879 & 0.0070\\
 0.500 & 0.193 & 5.194 & 0.689 & 0.853 & 0.836 & 0.0078\\
 0.500 & 0.193 & 5.194 & 0.933 & 0.866 & 0.862 & 0.008\\
 0.500 & 0.193 & 5.194 & 1.355 & 0.965 & 0.960 & 0.0078\\
 0.500 & 0.193 & 5.194 & 2.237 & 1.273 & 1.263 & 0.0075\\
 \hline
 0.500 & 0.289 & 7.354 & 0.179 & 2.489 & 2.438 & 0.0262\\
 0.500 & 0.289 & 7.354 & 0.400 & 2.329 & 2.251 & 0.0254 \\
 0.500 & 0.289 & 7.354 & 0.655 & 2.377 & 2.329 & 0.0268\\
 0.500 & 0.289 & 7.354 & 1.137 & 2.728 & 2.718 & 0.0256\\
 0.500 & 0.289 & 7.354 & 2.258 & 3.704 & 3.672 & 0.0212\\
 \hline
 0.500 & 0.347 & 11.746 & 0.231 & 6.365 & 6.403 & 0.1157\\
 0.500 & 0.347 & 11.746 & 0.558 & 6.624 & 6.349 & 0.1121\\
 0.500 & 0.347 & 11.746 & 1.227 & 8.089 & 7.813 & 0.0819\\
 0.500 & 0.347 & 11.746 & 2.882 & 11.288 & 11.189 & 0.0692\\
 \hline
 0.500 & 0.366 & 16.666 & 0.159 & 12.427 & 11.580 & 0.2998\\
 0.500 & 0.366 & 16.666 & 0.582 & 13.300 & 13.565 & 0.2851\\
 0.500 & 0.366 & 16.666 & 1.491 & 17.028 & 17.194 & 0.2082\\
 0.500 & 0.366 & 16.666 & 3.769 & 23.994 & 23.571 & 0.1409\\
\hline
\end{tabular}}
\end{table}

\maketitle
  \begin{table}[ht]
\caption{$N=4000$,  $\lambda$ dependence of the $\MSE$ at fixed $\mu$\label{table:FullVaryLambdaFourK}}  
\centering 
\vspace{.5cm}
\begin{tabular}{|l|l|l|l|l|l|l|l|l|l|l|l|} 
\hline 
$\delta$& $\rho$ & $\mu$ & $\lambda$ & fMSE & eMSE & SE\\
 [1ex] 
\hline 
0.100 & 0.095 & 5.791 & 0.402 & 0.152 & 0.144 & 0.0017\\
 0.100 & 0.095 & 5.791 & 1.258 & 0.136 & 0.128 & 0.0016 \\
 0.100 & 0.095 & 5.791 & 2.037 & 0.142 & 0.133 & 0.0016\\
 0.100 & 0.095 & 5.791 & 3.169 & 0.174 & 0.168 & 0.0016\\
 0.100 & 0.095 & 5.791 & 4.948 & 0.239 & 0.228 & 0.0012\\
 \hline
 0.100 & 0.142 & 8.242 & 0.804 & 0.380 & 0.348 & 0.0064 \\
 0.100 & 0.142 & 8.242 & 1.960 & 0.408 & 0.389 & 0.0058\\
 0.100 & 0.142 & 8.242 & 3.824 & 0.534 & 0.510 & 0.0051\\
 0.100 & 0.142 & 8.242 & 6.865 & 0.737 & 0.716 & 0.0034\\
 \hline
 0.100 & 0.170 & 12.906 & 0.465 & 1.045 & 0.950 & 0.0228\\
 0.100 & 0.170 & 12.906 & 2.298 & 1.178 & 1.111 & 0.0232\\
 0.100 & 0.170 & 12.906 & 5.461 & 1.619 & 1.591 & 0.0159\\
 0.100 & 0.170 & 12.906 & 10.607 & 2.197 & 2.182 & 0.008\\
 \hline
 0.100 & 0.180 & 18.278 & 0.338 & 2.063 & 1.588 & 0.0619\\
 0.100 & 0.180 & 18.278 & 2.934 & 2.467 & 2.171 & 0.0532\\
 0.100 & 0.180 & 18.278 & 7.545 & 3.474 & 3.367 & 0.0312\\
 0.100 & 0.180 & 18.278 & 14.997 & 4.677 & 4.551 & 0.0169\\
 \hline
 \hline
 0.150 & 0.109 & 5.631 & 0.420 & 0.236 & 0.228 & 0.0022\\
 0.150 & 0.109 & 5.631 & 1.073 & 0.212 & 0.209 & 0.0023\\
 0.150 & 0.109 & 5.631 & 1.700 & 0.218 & 0.213 & 0.0021\\
 0.150 & 0.109 & 5.631 & 2.657 & 0.260 & 0.251 & 0.0024\\
 0.150 & 0.109 & 5.631 & 4.284 & 0.359 & 0.353 & 0.0017\\
 \hline
  0.150 & 0.163 & 8.030 & 0.720 & 0.588 & 0.595 & 0.0072\\
 0.150 & 0.163 & 8.030 & 1.614 & 0.626 & 0.610 & 0.0078\\
 0.150 & 0.163 & 8.030 & 3.135 & 0.804 & 0.807 & 0.0058\\
 0.150 & 0.163 & 8.030 & 5.868 & 1.125 & 1.118 & 0.0047\\
 \hline
 0.150 & 0.196 & 12.577 & 0.434 & 1.612 & 1.572 & 0.0341\\
 0.150 & 0.196 & 12.577 & 1.814 & 1.792 & 1.720 & 0.0281\\
 0.150 & 0.196 & 12.577 & 4.339 & 2.433 & 2.383 & 0.0205\\
 0.150 & 0.196 & 12.577 & 8.903 & 3.359 & 3.333 & 0.0126\\
 \hline
 0.150 & 0.207 & 17.814 & 0.305 & 3.185 & 2.864 & 0.0861\\
 0.150 & 0.207 & 17.814 & 2.231 & 3.715 & 3.582 & 0.0722\\
 0.150 & 0.207 & 17.814 & 5.879 & 5.202 & 5.141 & 0.0439\\
 0.150 & 0.207 & 17.814 & 12.455 & 7.142 & 7.154 & 0.0269\\
\hline
\end{tabular}
\end{table}

\begin{table}
\caption{$N=1500$,  $\mu$ dependence of the $\MSE$ at fixed $\lambda$}
\begin{center}
\begin{tabular}{|l|l|l|l|l|l|l|}
\hline 
$\delta$& $\rho$ & $\mu$ & $\lambda$ & fMSE & eMSE & SE \\
 [1ex] 
\hline 
 0.100 & 0.095 & 5.291 & 1.253 & 0.131 & 0.125 & 0.0022\\
 0.100 & 0.095 & 5.541 & 1.256 & 0.134 & 0.132 & 0.0025\\
 0.100 & 0.095 & 5.691 & 1.257 & 0.135 & 0.126 & 0.0027\\
 0.100 & 0.095 & 5.791 & 1.258 & 0.136 & 0.129 & 0.0024\\
 0.100 & 0.095 & 5.891 & 1.259 & 0.137 & 0.125 & 0.0027\\
 0.100 & 0.095 & 6.041 & 1.260 & 0.138 & 0.126 & 0.0030\\
 0.100 & 0.095 & 6.291 & 1.262 & 0.139 & 0.127 & 0.0028\\
 0.100 & 0.095 & 6.791 & 1.264 & 0.141 & 0.125 & 0.0031\\
\hline
 0.100 & 0.142 & 7.242 & 0.794 & 0.349 & 0.317 & 0.0074\\
 0.100 & 0.142 & 7.742 & 0.800 & 0.366 & 0.335 & 0.0084\\
 0.100 & 0.142 & 7.992 & 0.802 & 0.373 & 0.351 & 0.0089\\
 0.100 & 0.142 & 8.000 & 0.802 & 0.373 & 0.362 & 0.0094\\
 \hline
 \hline
0.250 & 0.134 & 4.459 & 0.952 & 0.338 & 0.336 & 0.0036\\
 0.250 & 0.134 & 4.959 & 0.957 & 0.359 & 0.346 & 0.0040\\
 0.250 & 0.134 & 5.209 & 0.959 & 0.367 & 0.356 & 0.0044\\
 0.250 & 0.134 & 5.359 & 0.960 & 0.371 & 0.373 & 0.0049 \\
 0.250 & 0.134 & 5.459 & 0.961 & 0.374 & 0.362 & 0.0047 \\
 0.250 & 0.134 & 5.559 & 0.962 & 0.376 & 0.367 & 0.0045 \\
 0.250 & 0.134 & 5.709 & 0.962 & 0.379 & 0.372 & 0.0048\\
 0.250 & 0.134 & 5.959 & 0.963 & 0.383 & 0.362 & 0.0052\\
 0.250 & 0.134 & 6.459 & 0.964 & 0.387 & 0.387 & 0.0058\\
\hline
 0.250 & 0.201 & 6.683 & 0.587 & 0.939 & 0.899 & 0.0126 \\
 0.250 & 0.201 & 7.183 & 0.590 & 0.988 & 0.965 & 0.0147 \\
 0.250 & 0.201 & 7.433 & 0.591 & 1.009 & 0.956 & 0.0147 \\
 0.250 & 0.201 & 7.583 & 0.592 & 1.021 & 1.027 & 0.0155 \\
 \hline
 \hline
 0.500 & 0.193 & 4.194 & 0.684 & 0.769 & 0.770 & 0.0052\\
 0.500 & 0.193 & 4.694 & 0.687 & 0.818 & 0.823 & 0.0066\\
 0.500 & 0.193 & 4.944 & 0.688 & 0.837 & 0.838 & 0.0073\\
 0.500 & 0.193 & 5.094 & 0.689 & 0.847 & 0.835 & 0.0068\\
 0.500 & 0.193 & 5.194 & 0.689 & 0.853 & 0.834 & 0.0073\\
 0.500 & 0.193 & 5.294 & 0.689 & 0.858 & 0.845 & 0.0079\\
 0.500 & 0.193 & 5.444 & 0.690 & 0.865 & 0.863 & 0.0079\\
 0.500 & 0.193 & 5.694 & 0.690 & 0.874 & 0.887 & 0.0085\\
 0.500 & 0.193 & 6.194 & 0.691 & 0.886 & 0.868 & 0.0085\\
 \hline
 0.500 & 0.289 & 6.354 & 0.398 & 2.119 & 2.071 & 0.0195\\
 0.500 & 0.289 & 6.854 & 0.399 & 2.234 & 2.214 & 0.0235\\
 0.500 & 0.289 & 7.104 & 0.399 & 2.284 & 2.157 & 0.0252\\
 0.500 & 0.289 & 7.254 & 0.400 & 2.313 & 2.271 & 0.0244\\
 0.500 & 0.289 & 7.354 & 0.400 & 2.329 & 2.316 & 0.0275\\
 0.500 & 0.289 & 7.454 & 0.400 & 2.346 & 2.287 & 0.0287\\
 0.500 & 0.289 & 7.604 & 0.400 & 2.370 & 2.327 & 0.0306\\
 0.500 & 0.289 & 7.854 & 0.401 & 2.404 & 2.339 & 0.0284\\
 0.500 & 0.289 & 8.000 & 0.401 & 2.422 & 2.409 & 0.0300\\
\hline
\end{tabular}
\end{center}
\label{table:VaryMu}
\end{table}

  \begin{table}[ht]
\caption{$N=1500$, $\MSE$ for 5-point prior}\label{theomatch1500} 
\centering 
\vspace{.5cm}
\begin{tabular}{|l|l|l|l|l|l|l|l|l|l|l|} 
\hline 
$\delta$& $\rho$ & $\mu$ & $\lambda$ & Theoretical MSE & Empirical MSE & $\alpha$\\
 [1ex] 
\hline 
0.250  &0.134  &1.894  &0.857 &0.120  &0.151 & 0 \\
0.250  &0.134  &2.171  &0.897 &0.162  &0.163 & 0.122\\
0.250  &0.134  &2.447  &0.901 &0.178  &0.177 & 0.244\\
0.250  &0.134  &2.724  &0.906 &0.196  &0.195 & 0.366\\
0.250  &0.134  &3.001  &0.912 &0.215  &0.210 & 0.488\\
0.250  &0.134  &3.277  &0.918 &0.237  &0.236 & 0.611\\
0.250  &0.134  &3.554  &0.926 &0.261  &0.257 & 0.7333\\
0.250  &0.134  &3.830  &0.935 &0.287  &0.280 & 0.8556\\
0.250  &0.134  &4.107  &0.945 &0.317  &0.307 & 0.9778\\
0.250  &0.134  &4.383  &0.957 &0.348  &0.359 & 1.1000\\
\hline
\end{tabular}
\end{table}

\bibliographystyle{amsalpha}
\providecommand{\bysame}{\leavevmode\hbox to3em{\hrulefill}\thinspace}
\providecommand{\MR}{\relax\ifhmode\unskip\space\fi MR }
\providecommand{\MRhref}[2]{%
  \href{http://www.ams.org/mathscinet-getitem?mr=#1}{#2}
}
\providecommand{\href}[2]{#2}

\end{document}